\documentclass[12pt,smallextended,natbib,runningheads,epsfig]{article}
\usepackage{latexsym, amsfonts}
\usepackage{amssymb, amsmath}
\usepackage{amsfonts}
\usepackage{graphicx}
\usepackage{setspace} 
\usepackage{rotating}
\usepackage{mathtools}

\mathtoolsset{showonlyrefs}
 \usepackage{pdfsync}
\usepackage{tikz}
\usepackage{mathrsfs}
\usepackage{url}
\usepackage[authoryear]{natbib}
\newcommand{\be}{\begin{equation}}
\newcommand{\ee}{\end{equation}}
\newcommand{\bea}{\begin{eqnarray}}
\newcommand{\eea}{\end{eqnarray}}
\newcommand{\bean}{\begin{eqnarray*}}
\newcommand{\eean}{\end{eqnarray*}}
\newcommand{\brray}{\begin{array}}
\newcommand{\erray}{\end{array}}
\newcommand{\ben}{\begin{equation}{nonumber}}
\newcommand{\een}{\end{equation}{nonumber}}

\newtheorem{dfn}{Definition}[section]
\newtheorem{thm}[dfn]{Theorem}
\newtheorem{lmma}[dfn]{Lemma}
\newtheorem{ppsn}[dfn]{Proposition}
\newtheorem{crlre}[dfn]{Corollary}
\newtheorem{xmpl}[dfn]{Example}
\newtheorem{rmrk}[dfn]{Remark}
\newtheorem{algo}[dfn]{Algorithm}

\newcommand{\bdfn}{\begin{dfn}}
\newcommand{\bthm}{\begin{thm}}
\newcommand{\blmma}{\begin{lmma}}
\newcommand{\bppsn}{\begin{ppsn}}
\newcommand{\bcrlre}{\begin{crlre}}
\newcommand{\bxmpl}{\begin{xmpl}}
\newcommand{\brmrk}{\begin{rmrk}}
\newcommand{\edfn}{\end{dfn}}
\newcommand{\ethm}{\end{thm}}
\newcommand{\elmma}{\end{lmma}}
\newcommand{\eppsn}{\end{ppsn}}
\newcommand{\ecrlre}{\end{crlre}}
\newcommand{\exmpl}{\end{xmpl}}
\newcommand{\ermrk}{\end{rmrk}}

\newcommand{\bs}{\boldsymbol}

\newcommand{\E}{\mathbb{E}}
\newcommand{\HE}{\textcolor{blue}}

\newcommand{\bsf}{{\boldsymbol{\mathcal F}}}

\newtheorem{assumption}{Assumption}[section]
\newtheorem{theorem}{Theorem}[section]
\newtheorem{example}{Example}[section]
\newtheorem{remark}{Remark}[section]
\newtheorem{lemma}{Lemma}[section]

\newtheorem{corol}{Corollary}[section]

\newtheorem{proposition}{Proposition}[section]

\def\a*{{\cal A}_{h,*}}
\def\B{{\cal B}(h)}
\def\B1{{\cal B}_1(h)}
\def\b{{\cal B}^{\rm s.a.}(h)}
\def\b1{{\cal B}^{\rm s.a.}_1(h)}

\newcommand{\FF}{\mathcal{F}}

\newcommand{\lf}{\lfloor}
\newcommand{\rf}{\rfloor}

\begin{document}
  \title{Identifying  shifts between two regression curves}
\author{\small Holger Dette \\
\small Ruhr-Universit\"at Bochum \\
\small Fakult\"at f\"ur Mathematik \\
\small 44780 Bochum, Germany \\
\and
\small Subhra Sankar Dhar \\
\small IIT Kanpur\\
\small Department of Mathematics \& Statistics \\
\small  Kanpur 208016, India\\
\and
\small Weichi Wu \\
\small  Tsinghua University\\
\small  Center for Statistics\\
\small  Department of Industrial Engineering\\
\small 10084 Beijing China
}

\maketitle

\begin{abstract}
This article studies the problem whether two convex (concave) regression functions 
modelling the relation  between a response and covariate in two  samples 
differ by a shift in the horizontal and/or  vertical axis. We consider a  nonparametric situation
assuming only smoothness of the regression functions.
A graphical tool based on the derivatives of the regression functions and their inverses 
is proposed to answer this question and studied in several examples. 
We also formalize this question in a corresponding hypothesis
and develop  a statistical test.  The asymptotic properties  of the corresponding test statistic are
 investigated under the null hypothesis and  local alternatives. In contrast to most of the literature
 on comparing shape invariant models, which requires  independent data
 the procedure is applicable for dependent and   non-stationary 
 data.  We also illustrate the finite sample 
 properties of 
 the new test  by  means of a small  simulation study and  a real data example.      

\end{abstract}



\section{Introduction}
\label{sec1}
\def\theequation{1.\arabic{equation}}
\setcounter{equation}{0}

A common problem in statistical analysis is the comparison of   two regression models that relate a common response variable to the same covariates for two different groups.  If the two  regression functions coincide such  statistical  inference can be performed on the basis of the pooled sample
and therefore it is of interest to test hypotheses of this type.
 More  formally,  let
\begin{align}
Y_{i,1}  &= m_{1} (t_{i,1}) + e_{i,1 }~,~i=1, \ldots, n_{1}  \label{m1} \\
Y_{j,2} &= m_{2} (t_{j,2}) + e_{j,2}~,~j=1, \ldots, n_{2} \label{m2}
\end{align}
denote two regression models with real valued responses and  predictors  $t_{\ell,k}$ and  random errors  $e_{i,1}$ and $e_{j,2}$.  Statistical methodology addressing  the  question, if the two regression functions  $m_{1}$ and $m_{2}$ coincide,
has been investigated by many authors, and there exists an enormous amount of  literature 
addressing this  important  testing  problem [see, for example \cite{hallhart1990,dettmunk1998,dettneum2001,neumdett2003}
 for some early and \cite{vilarfernandez2007,neumpard2009,maity2012,degras2012,durgrolop2013,parkhannkang2014} for some more recent references among many others]. 

Another interesting question in this context is the comparison of  the regression curves up to a certain parametric
transformation. Such parametric relationship between two regression curves often can be fitted into various real life examples; for instance, as it is mentioned in \cite{haemar1990}, the growth curves of children may have a simple parametric relationship between them. It may happen that those curves are realizations of one curve but differ in the time and the vertical axes, and consequently, the difference among those set of regression curves can be measured by two unknown quantities, namely, the horizontal shift (i.e., along the covariate axis) and the vertical scale (i.e., along the response axis). 
  
Many authors  have  worked on this problem. Exemplary we mention the early work  by 
\cite{haemar1990,carrhall1993,ronn2001} and the more recent references \cite{gamloumaz2007,vimond2010,colldala2015} among others. 
Several  authors proposed  tests for the hypotheses that  the  regression curves coincide 
up to a certain parametric relationship. The  proposed methodology is  based on the estimation of the
parametric form from the given data. 
In this article we contribute to this literature and propose a simple method  to test the hypothesis 
\begin{align} \label{1.0}
H_{0} : m_{1} (x) = m_{2} (x + c) + d   ~~\text{ for some constants } c, d ~, 
\end{align}  
where $m_{1}$ and $m_{2}$ 
are convex (or concave)  functions. 
The assumption of a  convex or concave regression
function is well justified in several applications.
For example,  production functions  are  often  assumed  to  be  concave  [see \cite{varian1984}],
economic theory implies  that utility functions are  concave [see \cite{matzkin1991}] or in 
finance  theory  restricts  call  option prices to be convex [see \cite{sahaduar2003}].

We will show in Section \ref{sec2} 
 that under the null hypothesis    \eqref{1.0}  the functions $((m_1')^{-1})' $
and $((m_2')^{-1})' $ coincide (here and throughout this paper $f'$ denotes the derivative of
the function  $f$ and $f^{-1}$ its inverse). This fact  is utilized to develop a graphical device to check the assumption 
\eqref{Nullhypo} by estimating the difference $((m_1')^{-1})' - ((m_2')^{-1})' $.  For this purpose, 
we use ideas of  \cite{DNP2006} who proposed a very 
simple estimator of the inverse regression function  say $f$ based on a 
kernel density estimation of the random variable $f(U)$, where $U$ is  uniformly  distributed random variable
on the interval  $(0, 1)$, and $f$ is either $m_{1}'$ or $m_{2}'$.   

The second contribution of this paper is  a formal 
test for  the hypothesis \eqref{1.0}  in the context of dependent and non-stationary  data, which
is   based on a suitable distance between estimates  of the functions $((m_1')^{-1})'$ and $((m_2')^{-1}(t))'$. 
More precisely, we investigate an $L^{2}$-norm of
 a  smooth estimator of  the difference  $((m_1')^{-1})'-((m_1')^{-1})'$  
and derive the asymptotic distribution of the corresponding  test statistic under the null hypotheses and local alternatives.  
  The challenges in deriving these results are twofold. First  - in contrast to most of the literature - we allow for a  very 
   complex dependence structure of the errors in  models \eqref{m1} and \eqref{m2}.
   In particular they can be time dependent and non stationary 
   [see, for example \cite{DA1997}, \cite{mallat1998adaptive}, \cite{ombao2005slex},
   \cite{nason2000wavelet}, \cite{zhou2009local}, \cite{vogt2012nonparametric}  for various definitions of non-stationary time series].
A particular difficulty  consists in the proof of  the asymptotic distribution of the estimated  integrated squared difference, which is 
(after appropriate standardization) normal, but involves  higher  order derivatives of the regression functions. As these quantities are  very difficult to estimate   we develop  a bootstrap test, which has very good finite sample properties and  is 
based on a Gaussian approximation used in the proof of the weak convergence of the test statistic.

The rest of the article is organized as follows. Section \ref{sec2} describes the basic methodology adopted in this article. A 
 new graphical device is proposed  for comparing two non-parametric regression functions up a to shift in the  covariate and response
 in Section  \ref{sec21}. The formal testing problem is considered in Section \ref{sec22}, while we give some theoretical 
 justification for these tools in Section \ref{sec3}.
 A small simulation study is carried out in Section \ref{sec4}, illustrating the finite sample properties of the proposed method and an application is discussed in  Section 
 \ref{sec5}. Finally, all proofs except of the proof of Lemma \ref{lem1}, which justifies our approach, 
 are given in an appendix in Section \ref{sec6}.

\section{Methodology} 
\label{sec2}
\def\theequation{2.\arabic{equation}}
\setcounter{equation}{0}

Throughout this paper  we consider two data sets $\{Y_{i,1}\}_{i=1,..,n_1}$ and  $\{Y_{i,2}\}_{i=1,..,n_2}$ that can be modelled as
\begin{align}\label{Model_1}
    Y_{i,s} = m_{s} \Big (\frac{i}{n_s}\Big ) + e_{i,s},~~i = 1, \ldots, n_{s},~s = 1,2~, 
\end{align}  
 the error random variables $\{ e_{i,1} \}_{i = 1, \ldots, n_{1}}$ and $\{ e_{i,2} \}_{i = 1, \ldots, n_{2}}$ 
are 
locally stationary process  satisfying some technical conditions that will be described later in Section \ref{sec31}, and $m_{1}$  and $m_2$, are unknown  sufficiently smooth regression functions.
We assume that $m_{1}$ and $m_{2}$ are convex (the case of concave regression functions can be treated in a similar manner)
and  are interested to investigate  in a hypothesis 
\begin{equation}\label{Nullhypo}
~~~~~~~H_{0} :  
\left\{
\begin{split}
 & \text{there exists constants $c \in (0,1)$ and  $d\in\mathbb{R}$ such that}   \\ 
& m_{1} (t) = m_{2} (t + c) + d,~
\text{for all  $t\in (0,1-c)$ } 
\end{split}
\right.
~~~~~~
\end{equation}  
Notice that we assume that information about the sign of a potential vertical shift can be obtained by visible inspection of the data.  
A corresponding hypothesis with a vertical shift by a negative constant  $c$ can be formulated and treated in a similar way, but the details are omitted for the sake of brevity. 
 A key observation is that under the null hypothesis \eqref{Nullhypo} 
  we have \begin{align}\label{Keyobs}
  ((m_1')^{-1}(t))' - ((m_2')^{-1}(t))' = 0,
  \end{align} 
  and this fact motivates us to propose a test statistic and a graphical device based on the the estimate of $((m_1')^{-1}(t))' - ((m_{2}')^{-1}(t))'$.

 \begin{lemma} \label{lem1}
 Assume  that  the regression functions 
 $m_1$ and $m_2$ in \eqref{Model_1} 
 have a  strictly increasing first order derivative on the interval  $[0,1]$, then the following statements are equivalent. 
 \begin{description}
 \item (1)  There exists  a constant
$c\in (0,1)$ such that $m_1(t)=m_2(t+c)+d$ for all $t\in (0,1-c)$.
 \item (2)  Equation \eqref{Keyobs} holds for all  $u\in (m'_1(0), m'_1(1-c))$.
 \end{description}
 \end{lemma}
 
 {\it Proof.} If condition  (1) holds, then
 \begin{align}
     m_1'(t)=m_2'(t+c)
     \end{align}
     for all 
  $t \in (0,1-c)$.
Now consider the equation   $m_1' (x)=m_2'(x+c)=u$  for some fixed  $u\in(m'(0),m'(1-c))$ and   note that both derivatives are strictly increasing.
Consequently we obtain for a solution in the interval  for $(0,1-c)$
     \begin{align}
         x=(m_1')^{-1}(u)~;~~ x+c=(m_2')^{-1}(u)~.
     \end{align}
In particular, this yields (subtracting both equations)
\begin{align} \label{cdef}
c=(m_2')^{-1}(u)-(m_1')^{-1}(u)
\end{align}
for  any  $u\in(m_1'(0),m_1'(1-c))$. Taking derivatives on both sides of \eqref{cdef} gives \eqref{Keyobs} and 
shows that (1) implies (2).

On the other hand, if condition (2)  holds,
it follows
$$
        \int_{m'_1(0)}^s ((m_1')^{-1})'(u)du=\int_{m'_1(0)}^s ((m_2')^{-1})'(u)du,
$$
any $s\in (m'_1(0),m'_1(1-c))$,
which yields
$$
(m_2')^{-1}(s)=(m_1')^{-1}(s) +c ~
$$
   for $s\in (m_1'(0),m_1'(1-c))$, 
where 
  \begin{align}\label{value_c}c=( m_2')^{-1}( m'_1(0)).
  \end{align} 
Applying the function $m_2^\prime$ on both sides finally  
gives
$$
         m_2'((m_1')^{-1}(s)+c))=s=m_1'((m_1')^{-1}(s)) 
$$
 for $ s\in (m_1'(0),m_1'(1-c))$. 
 Using the notation $(m_1')^{-1}(s)=t$ and integrating  with respect to $t$  shows that this  is equivalent to (1), which completes the proof  of  Lemma \ref{lem1}.
 \hfill $\Box$

\subsection{Graphical Device}
\label{sec21}

According to Lemma \ref{lem1}, under null hypothesis, the points 
$$\{ (t,  f_{1} (t)-  f_{2} (t)) ~|~t \in(m_1'(0),m_1'(1-c))  \}$$ lie on 
the horizontal axis. 
In order to construct a graphical device, let $\hat f_{1}$ and $\hat f_{2}$ denote suitably chosen uniformly consistent  estimates of the functions   $f_1= (( m_1')^{-1}) ' $ and $f_2= (( m_2')^{-1} )' $, respectively,
let 
$\hat m_1^\prime  $ denote an estimate of the derivative  $m_1^\prime$, 
and let $\hat c $ be an estimate of the vertical 
shift $c$. We now   consider a collection of points 
\begin{equation} \label{cn}
{\cal{C}}_{n_1,n_2} =\{(t_{\ell},  \hat f_{1} (t_{\ell})- \hat f_{2} (t_{\ell}) )  \colon   t_{\ell} \in (\hat a + \eta, \hat b - \eta) ;~~ {\ell} = 1 , \ldots ,L \},
\end{equation} 
 where $\hat a = \hat{m_{1}'}(0) $ 
 and  $\hat{b} = \hat{m_{1}'}(1 - \hat{c}) $ are estimates of 
 ${m_{1}'}(0) $  and   ${m_{1}'}(1 - {c}) $, respectively, 
 $\eta$ is a small positive constant and   $L$ is a positive integer.  Under the null hypothesis,  
 the points  of  ${\cal{C}}_{n_1,n_2} $ should  cluster around the horizontal axis.

 Here the necessary estimates  can be constructed in various ways.  For example,  $\hat f_{1}$ and $\hat f_{2}$ can be obtained 
using a smooth nonparametric estimate of the derivative of the regression function
 and calculating the derivative of its inverse. 
The  inversion of the  nonparametric estimates  of the derivatives $m_{1}$ and $m_{2}$ might be difficult as these 
functions are usually not monotone. Possible solutions are to construct isotone (smooth) nonparametric estimates 
 of the derivatives as proposed in
  \cite{mammen1991} and \cite{hallhuang2001} among others and then calculate the inverse. 
 Here we  use a  more direct  approach   related to the work of  \cite{DNP2006}
 who   proposed methodology for nonparametric estimation of a  monotone  regression 
 function based on monotone rearrangements. 

 To be precise, let $K$ denote a kernel function,
 $ b_{n,1}$, $b_{n,2}$ two bandwidths and 
define the estimate of the regression function $m_s$ and its derivative $m_s^\prime$  for $t\in [b_{n,s},1-b_{n,s}]$  by  
\begin{align}\label{2019-11}
(\hat m_s(t), b_{n,s}\hat m'_{s}(t))^\top=\underset{\beta_0,\beta_1}{\mathrm{argmin}} \sum_{i=1}^n\Big(Y_{i,s}-\beta_0-\beta_1\Big(\frac{i}{n_s}-t\Big)\Big)^2K\Big(\frac{i/n_s-t}{b_{n,s}}\Big)~,~~ s=1,2,
\end{align}
and $\hat m'_s(t)=\hat m'_s(b_{n,s})$ for $0\leq t\leq b_{n,s}$, while  $\hat m'_s(t)=\hat m'_s(1-b_{n,s})$ for $1-b_{n,s}\leq t\leq 1.$ Let $K_d$ be  a kernel function,
 $h_{d}$ a sufficiently small bandwidth and $N$ a large positive integer (note that this is not the sample size).
We define the estimates
\begin{align} \label{f1} 
\hat f_{1}(t) & = \frac{1}{Nh_{d,1}}\sum\limits_{i = 1}^{N}K_d\Big (\frac{\hat{m}_{1}'(\frac{i}{N}) - t}{h_{d,1}}\Big ),  \\
 \label{f2} 
\hat f_{2}(t) & = \frac{1}{Nh_{d,2}}\sum\limits_{i = 1}^{N}K_d\Big (\frac{\hat{m}_{2}'(\frac{i}{N}) - t}{h_{d,2}}\Big ).
 \end{align} 
 for $f_1(t) = {(( m_{1}')^{-1})'} (t)$ and $ f_2(t)= {(( m_{2}')^{-1})'} (t)$, respectively.
 For the motivation of this definition note that, if the estimates $\hat{m}_{s}'$ 
 are  consistent for ${m}_{s}'$ ($s=1,2$), then we  can replace for  a sufficiently  large sample size the estimates by the unknown
  regression functions, and obtain by  
 a Riemann approximation
  (if $N \to \infty$, $h_{d}\to 0$)   
  \begin{align} 
\hat f_{s }(t)  &\approx  \frac{1}{Nh_{d}}\sum\limits_{i = 1}^{N}K_d\Big (\frac{{m}_{s}'(\frac{i}{N}) - t}{h_{d}}\Big )
\approx  \frac{1}{h_{d}} \int_{0}^{1}K_d\Big (\frac{{m}_{s}'(x) - t}{h_{d}}\Big ) dx  \\
  &=  \int_{({m}_{s}'(0) - t))/{h_{d}}}^{({m}_{s}'(1) - t))/{h_{d}}}  K_{d}(u) (({m}_{s}')^{-1})'(t + uh_{d}) du\approx  (({m}_{s}')^{-1})'(t)  \mathbf{1} \{ {m}_{s}'(0) <  t  < {m}_{s}'(1) \}. 
  \nonumber  
  \end{align} 
  where $ \mathbf{1}  (A)$ denotes the indicator functions of the set $A$ and we have used the fact that  ${m}_\ell'$  
  is   non-decreasing 
  (see  \cite{DNP2006} for more details). 
  Finally, 
  the estimate of $(m_2^\prime)^{-1}$ can be obtained by   integration, that is 
\begin{equation}
     \label{h15} 
     \hat g_2(x) = \int^x_{m_2'(0)}  \hat f_2 (t) dt  
  \end{equation}
  and using \eqref{cdef} we obtain an estimate 
    \begin{align}     \label{chat} 
  \hat c= \frac{1}{1-\tilde c}\int_0^{(1-\tilde c)}(\hat g_2(\hat m_1'(u))-u)du. 
  \end{align}
  of the vertical shift $c$.  
Here  $  \hat m_1^\prime $ is the estimate of the derivative of $m_1$ defined in \eqref{2019-11} and  $$
\tilde c=\hat g_2 (\hat m_1'(0)).
$$ 
is a preliminary  consistent estimator of $c$.
 The resulting estimates for $ a=m_1'(0)$  and   $b=m_1(1-c) $
are then given by 
\begin{eqnarray} 
\hat a&=&\hat m_1'(0)~,~~
\label{aest} 
\hat b = \hat m_1'(1-\hat c)
\end{eqnarray}
(note that we assume that $c>0$).
We will prove in Theorem \ref{The0} below
 that under  the null hypothesis  \eqref{Nullhypo}  the points of the set ${\cal{C}}_{n_1,n_2}$ will concentrate 
 around the horizontal axis when the sample sizes are sufficiently large.
Therefore we propose  a graphical device that plots the points of the set ${\cal{C}}_{n_1,n_2}$.

\begin{example} \label{ex1} 
{\rm
{
 We consider the  regression models  \eqref{Model_1} with independent standard normal distributed errors and different regression functions\, where the sample sizes are  $n_1 = n_2 = 100$. In this numerical study, $N = 100$, $h_{d, N} = N^{-1/3}$, and bandwidths $b_{n_1,1}$ and  $b_{n_2,2}$ are  chosen as described 
in Section \ref{sec4}. The set $\mathcal{C}_{n_1,n_2}$ consists of $L = 1000$ equally spaced  points  from the interval $(\hat a + \eta, \hat b - \eta)$, where  $\eta = 0.01$. To compute the 
local linear estimators  we use the   R package named `locpol'. The following models are considered in  this example: }


\begin{align}
\label{ex1h0} m_{1} (x) &=  (x - 0.4)^{2} \text{ ~~and  ~} m_{2} (x) = (x - 0.3)^{2} - 0.2 ,  \\
\label{ex1h1} m_{1} (x) &=  (x - 0.4)^{2} \text{ ~~and  ~} m_{2} (x) = x^{3} ,  \\
\label{ex3h0}
 m_{1}(x) &=  \sin (-\pi x)  \text{ ~~and  ~ } m_{2}(x) =\sin(-\pi(x + 0.1)) + \frac{1}{4}, \\
\label{ex3h1}
 m_{1}(x) &=  \sin (-\pi x)  \text{ ~~and  ~ } m_{2}(x) = -\cos (\pi x).
 \end{align}

Note that examples \eqref{ex1h0} and \eqref{ex3h0} 
correspond to the null hypothesis, while  \eqref{ex1h1} and \eqref{ex3h1} represent alternatives. 
The corresponding plots of the set $\mathcal{C}_{n_1,n_2}$  are shown  in Figure \ref{Fig1}, where the the left  panels  clearly support 
the null hypothesis   of a vertical and horizontal shift  between the regression functions (the points  are clustered around the $x$-axis).
On the other hand, the panels  on the right give clear evidence that  the null hypothesis \eqref{Nullhypo} is not true.

\begin{figure}\label{Fig1}
\includegraphics[width= 3in, height = 2.5in]{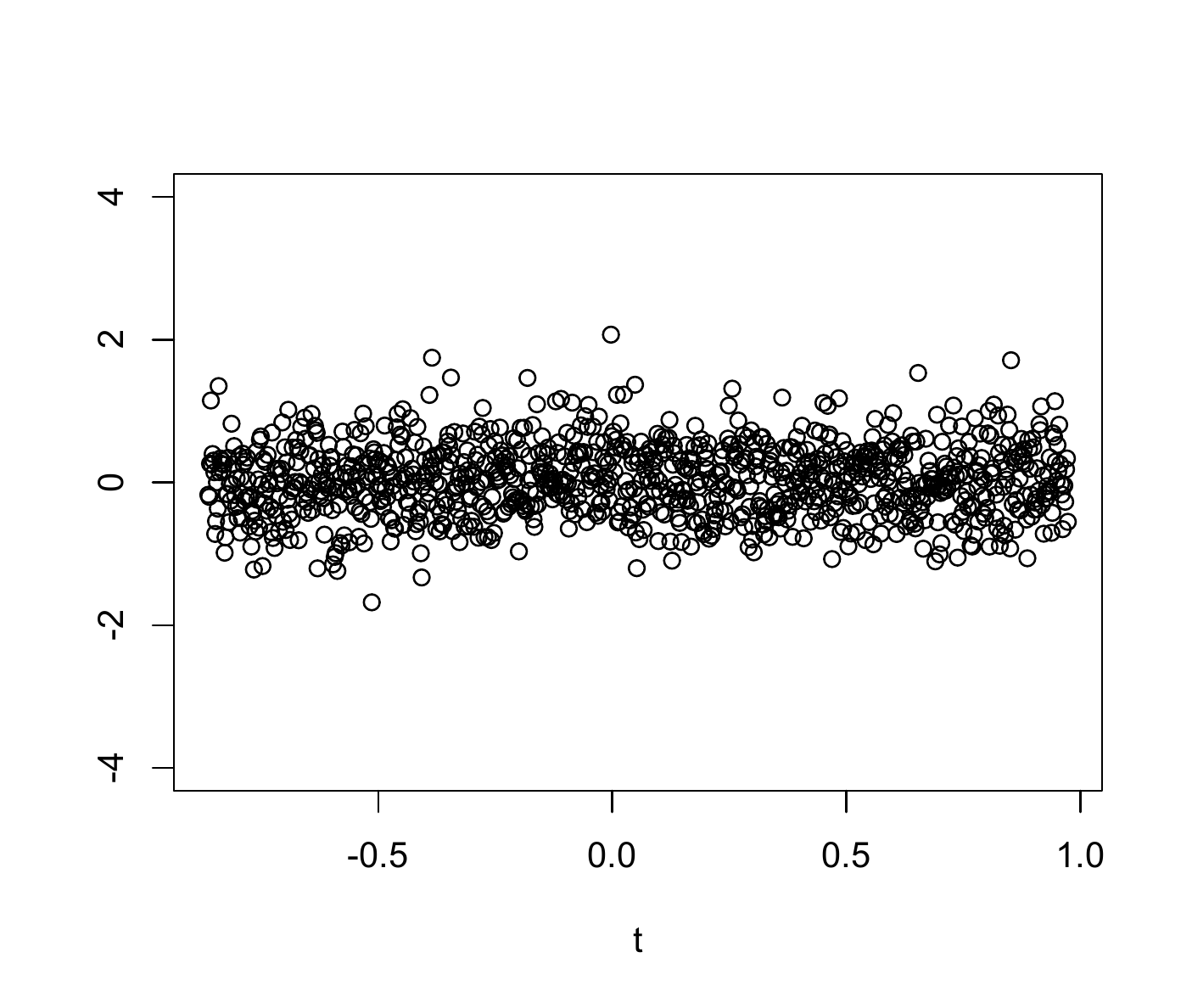}\includegraphics[width= 3in, height = 2.5in]{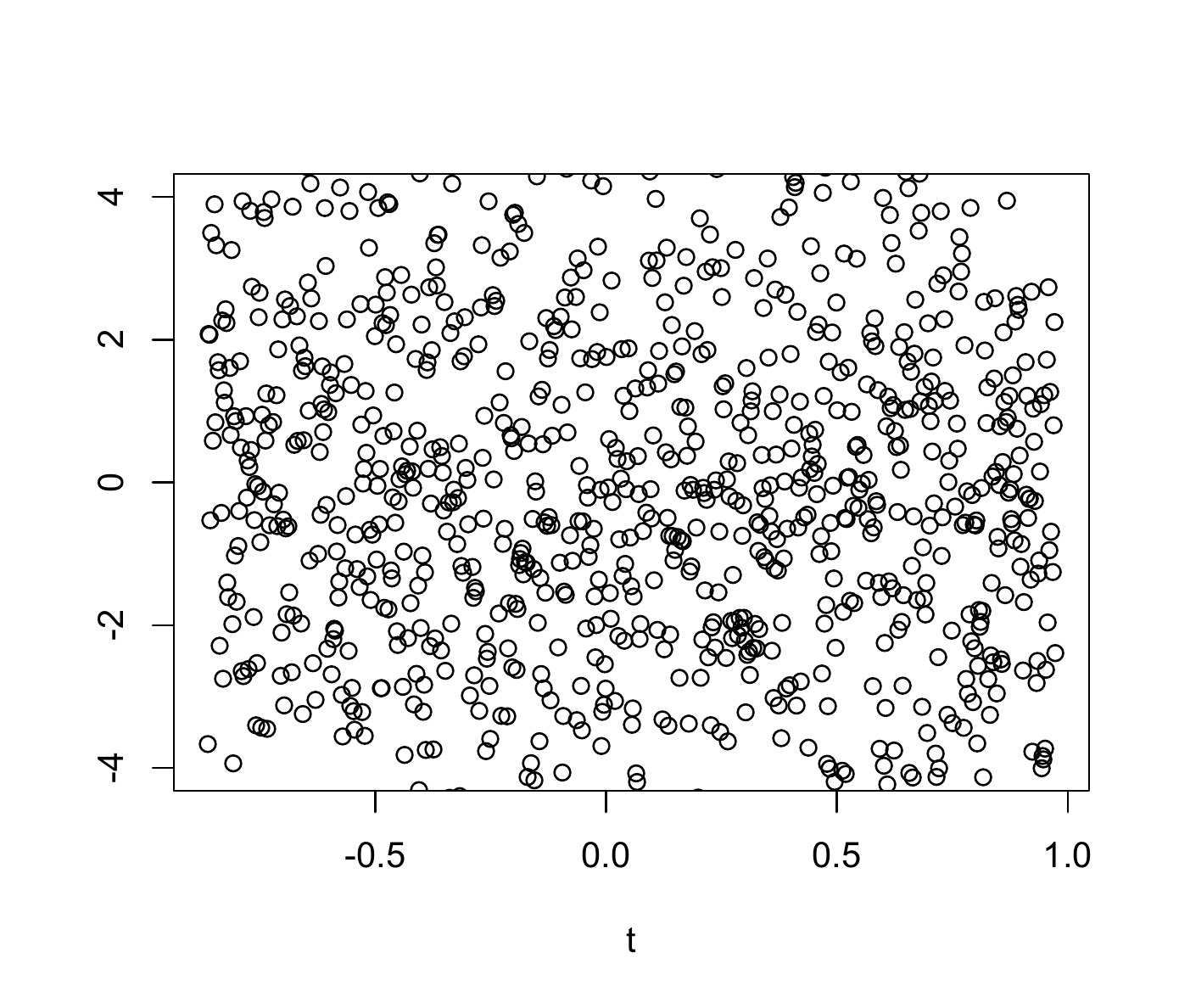}\\ 
\includegraphics[width= 3in, height = 2.5in]{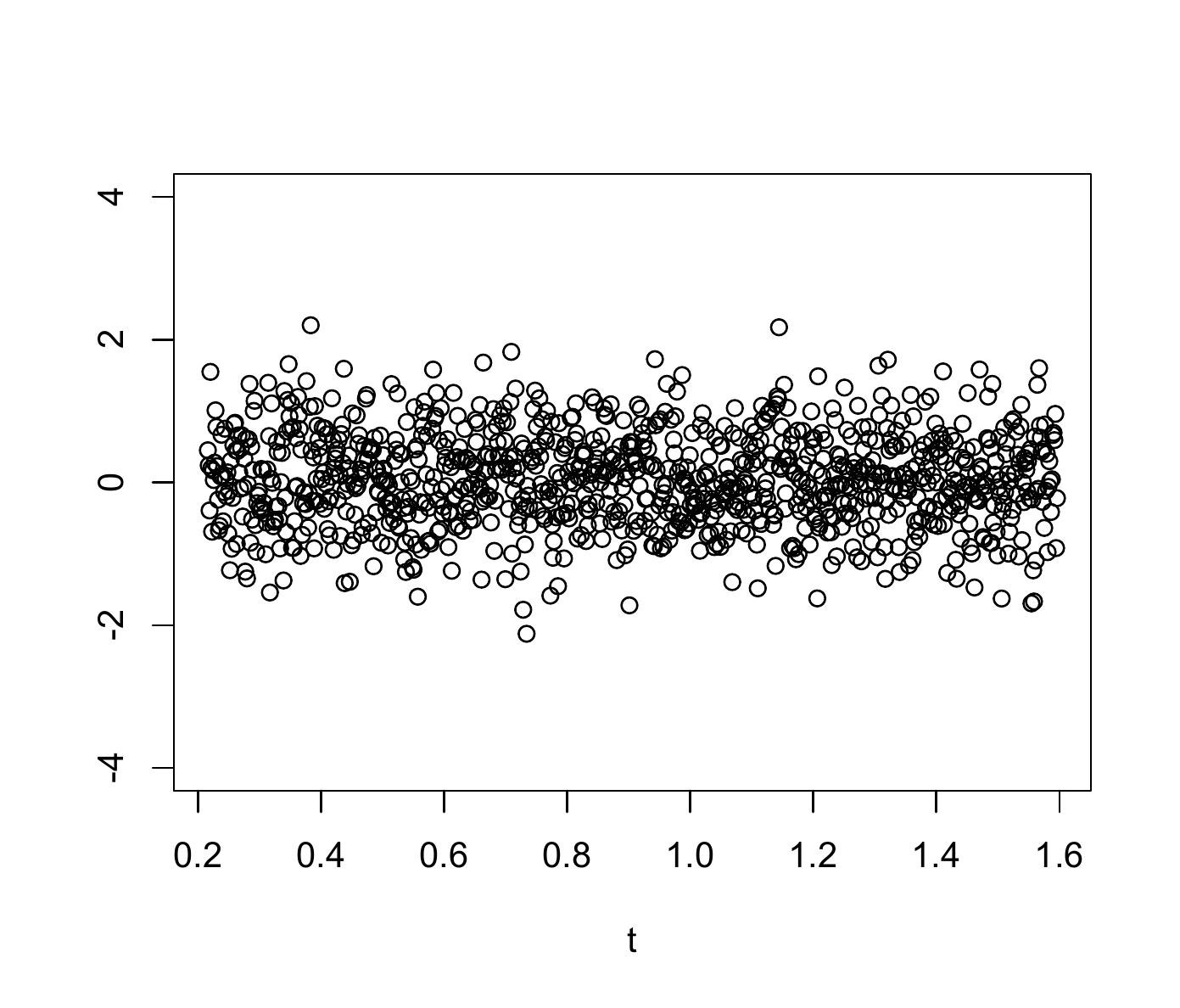}\includegraphics[width= 3in, height = 2.5in]{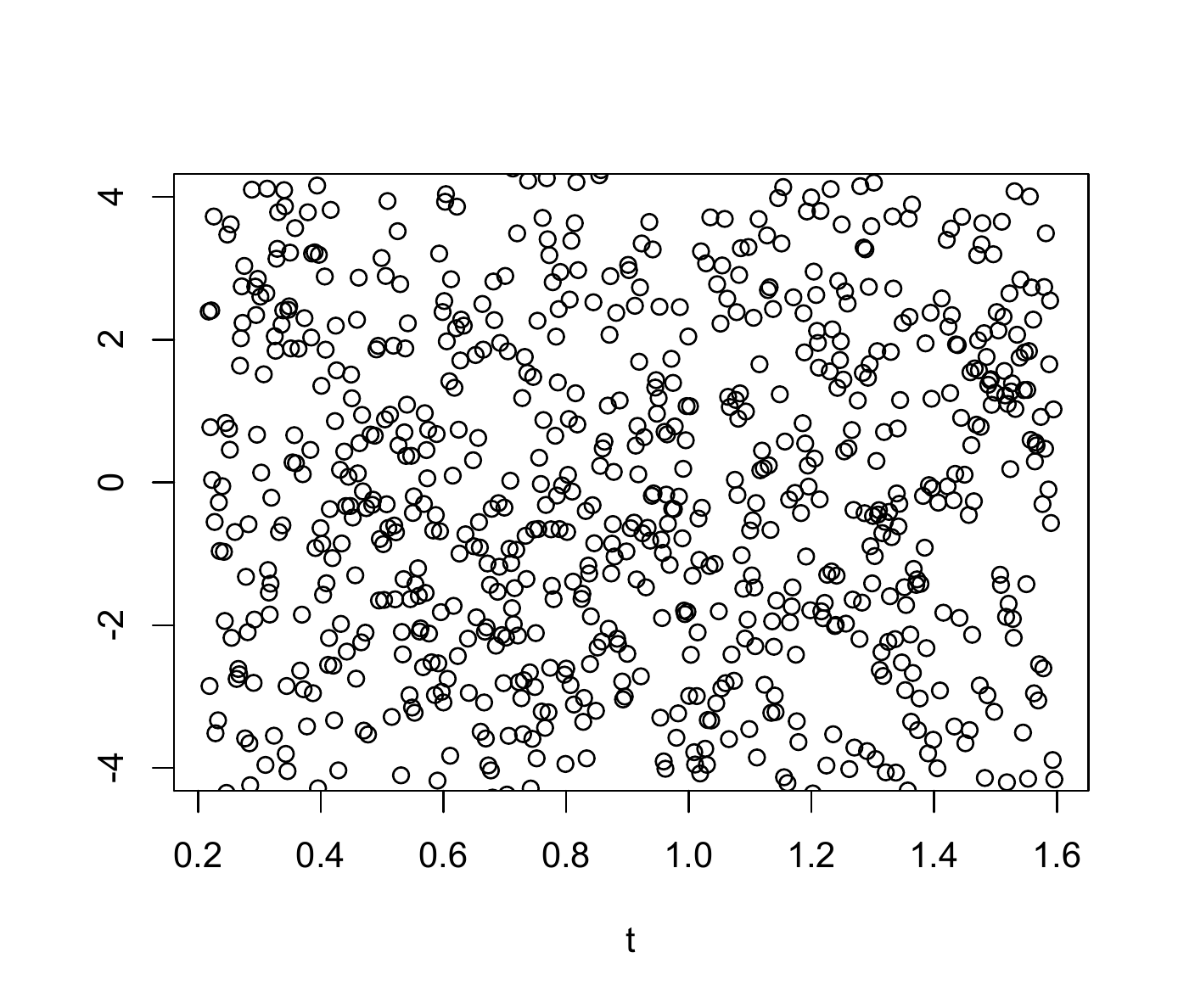}\\
\caption{
\it  
 Plots of the set ${\cal C}_{n_1,n_2} $ for different examples. The panels on the left correspond to the models \eqref{ex1h0}  and
 \eqref{ex3h0} (null hypothesis) and the panels on the right correspond to the models \eqref{ex1h1} and
 \eqref{ex3h1} (alternative).
 } 
\end{figure}
}
\end{example}

\subsection{Investigating shifts in the regression functions by testing }
\label{sec22}
The graphical device discussed in  the previous section provides a simple tool of visual examination of the null hypothesis \eqref{Nullhypo}, but does 
not give any information about the statistical uncertainty of a decision. In this section we 
will add to this tool a statistic which can be used
to rigorously test the null hypothesis \eqref{Nullhypo} at a controlled type I error. 
Recalling the definition of the estimates \eqref{f1} and \eqref{f2} of  ${( ( m_{1}')^{-1})'} (t)$ and ${( (  m_{2}')^{-1})'} (t)$, we propose to  reject the null  hypothesis \eqref{Nullhypo} for large values of 
the statistic 
\begin{align}\label{def_T}
T_{n_1,n_2}=\int\big(\hat f_1(t)-\hat f_2(t)\big)^2\hat w(t)dt,
\end{align}
where the weight function is defined by 
  \begin{align}
\hat w(t)=\mathbf 1(\hat a+\eta\leq t\leq \hat b-\eta),
\end{align}
$\eta $ is a small positive constant  
and   $\hat a$  and   $\hat b$ are defined in \eqref{aest}.
In fact, $\hat w(t)$ is a consistent estimator of 
the deterministic weight function 
\begin{align}\label{aandb}
w(t)=\mathbf 1(a+\eta\leq t\leq b-\eta),
\end{align}
where $  a=m_1'(0), b=m_1'(1-c). $

 \begin{remark} {\rm 
For the  construction of  the test statistic, 
other distances 
between  the functions  ${( ( \hat m_{1}')^{-1})'} (t)$ and ${( ( \hat m_{2}')^{-1})'} (t)$ could be considered as well. For the $L^2$ distance,  
the derivation of the asymptotic distribution of the statistic $T_{n_1,n_2}$ is already very  complicated (see Section \ref{sec6} for details), but we can make use
of a central limit theorem for random quadratic forms [see \cite{de1987central}]. Other distances such as the supremum or $L^1$ distance could be considered as well with additional technical arguments. 
}
 \end{remark}

\section{Asymptotic properties} \label{Section3}
\label{sec3}
\def\theequation{3.\arabic{equation}}
\setcounter{equation}{0}

 Before stating the asymptotic distribution of $T_{n_1,n_2}$, a few concepts and assumptions are stated for model \eqref{Model_1}. 
For the dependence structure, we use a common concept non-stationarity, which will be described first.

\subsection{Locally stationary processes and  basic assumptions } \label{sec31}
Recall the definition of  model \eqref{Model_1} and denote by  
$\{\mathbf e_i\}_{i\in \mathbb N} = 
\{ (\varepsilon_{i,1},\varepsilon_{i,2})^\top\}_{i\in \mathbb N} $
the vector of errors.  Note that $\{ 
\mathbf e_i\}_{i \in \mathbb N} $ defines a triangular array although this is not reflected in our notation. In particular we 
 assume $\{\mathbf e_i\}_{i \in \mathbb N}$ is a locally stationary process in the sense of \cite{zhou2009local}
 such that it has the form
\begin{align} \label{loc}
\mathbf e_i=\mathbf G(i/n,\bs{\FF}_i)=((G_1(i/n,\FF_i),G_2(i/n,\mathcal {G}_i)^\top), 1\leq i\leq n,
\end{align}
where $\bs G:[0,1]\times \mathbb R^{\infty}\rightarrow \mathbb R^2$ is a measurable nonlinear filter, 
$\bs{\FF}_i=(...,\boldsymbol \epsilon_{i-1},\bs\epsilon_i)$ is a filtration 
and  $\{\bs \epsilon_i=(\varepsilon_{i,1},\varepsilon_{i,2})^\top\}_{i\in \mathbb N} $ a sequence of independent identically distributed random variables. 
In \eqref{loc} $G_1$ and $G_2$ are the marginal filters  and  $\FF_i=(....,\varepsilon_{i-1,1},\varepsilon_{i,1})$, $\mathcal {G}_i=(...,\varepsilon_{i-1,2},\varepsilon_{i,2})$.
Moreover for any $p$-dimensional vector $\mathbf v=(v_1,...,v_p)^\top$ we define $|\mathbf v|=\sqrt{\sum_{i=1}^pv_i^2}$,  $\| \mathbf v \|_4 = (\E(|\mathbf v|^4))^{1/4}$
and
make the following basic assumptions.

\begin{assumption}\label{Assum3.1}\ 
	\begin{description}
		\item(a) $\mathbb E(\mathbf G(t,\bsf_0))=0$ ~\mbox{for} $t\in[0,1]$,~\mbox{and}  $\displaystyle\sup_{t\in[0,1]}\|\mathbf G(t,\bsf_0)\|_4<\infty$.
		\item(b) $\displaystyle\sup_{0\leq t<s\leq 1}\|\mathbf G(t,\bsf_0)-\mathbf G(s,\bsf_0)\|_4<\infty$.
		\item(c) Let  $\{\bs \epsilon_i^*\}_{i \in \mathbb N} $  denote an  independent  copy of 
		$\{\bs \epsilon_i\}_{i \in \mathbb N} $ and  define 
		the  filtration $\bsf_i^*=(\bs \epsilon_{-\infty},...,\bs \epsilon_{-1},\bs \epsilon_0^*,...,\bs \epsilon_i)$.
		There exists a constant $\rho\in(0,1)$ such that for any $k\geq 0$,
		\begin{align*}
		\delta_4(k):=\sup_{t\in[0,1]}\|\mathbf G(t,\bsf_k)-\mathbf G(t,\bsf_k^*)\|_4=O(\rho^k) ~.
		\end{align*}
			\item(d) There exists a constant $\nu_0>0$ such that the $2\times 2$ matrix $\Sigma^2(t)-\nu_0 I_2$ is strictly positive definite  for any $t\in[0,1]$, where $I_2$  is the $2 \times 2 $ identity matrix,  and $\Sigma^2(t)$ is the long run variance of the locally stationary process defined as
			\begin{align}
				\Sigma^2(t)=\sum_{s=0}^\infty \mathbb E\left(\mathbf{G}(t,\bs{\FF}_0)\mathbf G(t,\bs{\FF}_s)^\top\right).
			\end{align}
			\item (e) $\Sigma^2(t)$ is a diagonal matrix with  entities  $\sigma_1^2(t)$ and  $\sigma_2^2(t)$ (the long-run variances of process $G_1(\cdot,\FF_i)$ and 
			$G_2(\cdot,\mathcal G_i)$). 
	\end{description}
	\end{assumption}
	Note that  it follows from the definition   of $\delta_{4}(k)$ that  $\delta_4(k)=0$ for $k\leq 0$.
Assumptions (d) and (e) ensure that 
$\sigma^2_1(t)$  and $\sigma^2_2(t)$    are  non-degenerate such that $\displaystyle\inf_{t\in[0,1]}\sigma_s^2(t)>0$ ($s=1,2$).

Recalling the definition of the  
local linear estimator for the  derivatives $m_1'$ and $m_2'$ in  \eqref{2019-11} we make the following assumptions.

\begin{assumption}\label{AssumptionK}
~
{\rm
	\begin{description}
    \item (a) The kernel $K$ is  a symmetric and twice differentiable function with  compact support, say  $[-1,1]$. Furthermore, $\int_{-1}^1 K(x)dx=1$。 
		
		\item (b) The kernel $K_d$ is an even density  with  compact support, say  $[-1,1]$. 
		\end{description}
		}
		
	\end{assumption}
	
\begin{assumption}\label{AssumptionMean}
~
{\rm
	\begin{description}
 
 \item (a) 
		$m_1, m_2\in \mathcal C^{2,1}[0,1]$,  where 
		$\mathcal C^{2,1}[0,1]$ represents the  set of twice continuously differentiable functions, whose second order derivative is Lipschitz continuous on the interval $[0,1]$.
		\end{description}
		}
		
	\end{assumption}

\begin{assumption}\label{band}
{\rm
For $s=1,2$ let  
$$
\pi_{n,s}=\frac{\log n}{\sqrt{nb_{n,s}}b_{n,s}}+\frac{n^{1/4}\log ^2n}{nb_{n,s}^2}+b_{n,s}^2 ~,~~ \pi_{n,s}'=\frac{n^{1/4}\log ^2n}{nb_{n,s}^2}+b_{n,s}^2$$ 
and assume that $\pi_{n,s}=o(h_{d, n})$ ($s=1,2$). Further, assume that 
\begin{eqnarray*}
&&nb_{n,s}^2\rightarrow \infty~, ~~nb_{n,s}^{4}\log n \Big (\frac{\pi_{n,s}'}{b_{n,s}}+\frac{\pi_{n,s}^3}{h_d^3}+h_d+\frac{1}{Nh_d} \Big )^2=o(1), \\
&&
\bar \omega_n b_{n,s}^{-1/2}\log^2n =o(1)~,
\end{eqnarray*}
where 
\begin{align} \label{omega}
\bar \omega_{n,s}=\frac{\log n}{\sqrt{nb_{n,s}}b_{n,s}}+\frac{n^{1/4}\log ^2n}{nb_{n,s}^2}+b_{n,s}~,~~s=1,2.
\end{align}
}
\end{assumption}

 \subsection{Asymptotic properties of $\mathcal{C}_{n_1,n_2}$}
 \label{sec33}

The following theorem  describes the asymptotic properties of  the set  ${\cal{C}}_{n_1,n_2}$
defined in \eqref{cn} if it is used with the {local linear}  estimates \eqref{2019-11} for the derivatives 
$m_1'$ and $m_2'$. It basically gives a theoretical justification for the use of the graphical device
proposed in Section~\ref{sec21}.
The proof can be found in Section \ref{sec61a}.

\begin{theorem}\label{The0}
 Define for  $\epsilon > 0$ the set 
 $$L(\epsilon,g) = \{(x, y) : x\in [m_1'(0)+\eta, m_1'(1-c)-\eta], |y-g(x)|\leq\epsilon\}.
 $$
 where  $ g= ((m'_1)^{-1})' -((m'_2)^{-1})' $.  
If  Assumptions \ref{Assum3.1}--\ref{band}  are satisfied,
then we have 
$$
\lim_{n_1,n_2\rightarrow\infty}\mathbb P[{\cal{C}}_{n_1,n_2}\subset L(\epsilon,g)]= 1.
$$ 
\end{theorem}

\noindent
Under the null hypothesis we have $g \equiv 0 $ and 
$$L(\epsilon)  :=  L(\epsilon, 0)  = \{(x, y) : x\in [m_1'(0)+\eta, m_1'(1-c)-\eta], |y|\leq\epsilon\}.
$$
Theorem \ref{The0} shows,  that for large sample size  the  
points in the set ${\cal{C}}_{n_1,n_2}$ cluster around the horizontal axis if and only if  the null hypothesis holds.

\subsection{Weak convergence of the test statistic}
In this section,  we derive  the asymptotic distribution of the statistic $T_{n_1,n_2}$. For this purpose,  we  define
\begin{align} \label{kequiv}
    K^\circ(x)=\frac{K(x)x}{\int_{-1}^1K(x)x^2dx}, 
\end{align}
and obtain the following result. 
The proof is complicated and can be found in Section \ref{sec62}.

\begin{theorem} \label{thm32}
 Suppose that Assumption \ref{Assum3.1}-\ref{band} hold,  $n_2/n_1\rightarrow c_2$ for some constant $c_2 \in (0,\infty) $ and assume 
additionally that 
$$
\frac{b_{n,1}}{b_{n,2}}\rightarrow r_2\in (0,\infty).
$$
Consider local  alternatives of the form
$$
((m'_1)^{-1})'(t)-((m'_2)^{-1})'(t)=\rho_n g(t)+o(\rho_n),
$$
where 
$g \in \mathcal C[a,b]$,  $\rho_n=(n_1b_{n,1}^{9/2})^{-1/2}$ and the order  $o(\rho_n)$ of the remainder holds uniformly with respect to $t$. 
Then as $n_1,n_2\rightarrow \infty$,
\begin{align}\label{T_{n_1,n_2}asymptotics}
n_1b_{n,1}^{9/2}T_{n_1,n_2}- B_n(g) \Rightarrow {\cal N} (0, V_T),
\end{align} 
where the asymptotic bias and variance are given by
\begin{align*}
B_n(g) & = \frac{(\int_{-1}^1 vK'_d(v)dv)^2}{\sqrt{b_{n,1}}}((K^\circ)'*(K^\circ)'(0))\sum_{s=1}^2c_sr_s^{5}\int_{\mathbb R} \sigma^2_s(u)w(m_s'(u))(m_s''(u))^{-3}du  \\
& ~~~~~~~~~~~~~~~~~~~~~~~~~~~~~~~~~~~~~~~~~~~~~~~~~~~~~~~~~~   + \int_0^1g^2(t)w(t)dt\notag , \\
V_T &= 2 \Big (\int_{-1}^1 vK'_d(v)dv \Big )^4\sum_{s=1}^2c_s^2r_s^9\int_{\mathbb R} ((K^\circ)'*(K^\circ)'(z))^2dz
 \int_{\mathbb R}(\sigma_s^2(u)w(m_s'(u))(m''_s(u))^{-3})^2 du，\end{align*} 
$c_1=1$, $r_1=1$ respectively,  and $(K^\circ)'*(K^\circ)'$ denotes the convolution of the functions $(K^\circ)'$ and $(K^\circ)'$.
\end{theorem}
 
 \begin{remark} \label{rem31}
 { \rm Under the null hypothesis, we have $g \equiv0$
 and Theorem \ref{thm32}  can be  used to construct a 
 consistent asymptotic level $\alpha$ test for the hypotheses in \eqref{Nullhypo}. More precisely, the null hypothesis is rejected whenever 
 $$
 T_{n_1,n_2}>\frac{\hat B_n(0)+z_{1-\alpha} \hat V^{\frac{1}{2}}_T}{n_1b_{n,1}^{\frac{9}{2}}}~,
 $$
 where $z_{1-\alpha}$ is the  corresponding $(1-\alpha)$-th quantile, and 
 $\hat B_n(0)$ and $\hat V_T$ are  appropriate estimates of the asymptotic  bias (for $g(t) \equiv 0$) and variance, respectively. 
 Moreover, Theorem \ref{thm32} also  shows that this
 test is able to detect  alternatives converging to the null hypothesis at a  rate $\rho_n=(n_1b_{n,1}^{9/2})^{1/2}$. In this case,  
the asymptotic power of the test  is approximately given by  
  	\begin{align}
  	\Phi\Big ( \frac{\int g^2(t)w(t)dt}{V_T^{1/2}}-z_{1-\alpha}\Big ), 
  	\end{align}
  	where $\Phi $ is the cumulative distribution function of the standard normal distribution,
  	}
\end{remark}

In the case where the sample sizes $n_1$ and $n_2$ are equal Theorem \ref{thm32} directly leads to the following corollary.

 \begin{corol}
 \label{Thm1}   If the assumptions  of Theorem \ref{thm32} are satisfied,  the sample sizes  and bandwidths are equal (i.e. $n_1=n_2$ $b_{n,1}=b_{n,2}=b_{n}$ ), the weak convergence  in \eqref{T_{n_1,n_2}asymptotics} holds with
 \begin{align*}
     B_n(g) & =\frac{(\int vK'_d(v)dv)^2}{\sqrt{b_n}}((K^\circ)'*(K^\circ)')(0)\sum_{s=1}^2\int_{\mathbb R} \sigma^2_s(u)w(m_s'(u))(m_s''(u))^{-3}du \\
     & ~~~~~~~~~~~~
     ~~~~~~~~~~~~~~~
     ~~~~~~~~~~~~~~~~
     -\int_0^1g^2(t)w(t)dt\\
     V_T& =2\Big ( \int_{-1}^1 vK'_d(v)dv \Big)^4\sum_{s=1}^2\int_{\mathbb R} ((K^\circ)'*(K^\circ)'(z))^2dz\int_{\mathbb R}
     (\sigma_s^2(u)w(m_s'(u))(m''_s(u))^{-3})^2 du.
 \end{align*}
\end{corol}

\section{Implementation and simulation study}
\label{sec4}
\def\theequation{4.\arabic{equation}}
\setcounter{equation}{0}

We begin with some details regarding the implementation of the test. 
The calculation  of the test statistic requires 
the specification of the bandwidths and
we  use the  general  Cross Validation (GCV) method proposed in \cite{zhou2010simultaneous}. Specifically, let $\hat m_s(\cdot,b)$ denote the estimate of the regression function $m_s$ with bandwidth $b$, then we consider  
\begin{align*}
    \hat b_{n_s,s}=\text{argmin}_b \frac{n_s^{-1}\sum_{i=1}^{n_s}(Y_{i,s}-\hat m_s(i/n_s,b))^2}{(1-K(0)(n_sb)^{-1})^2}. 
\end{align*}
As pointed out by \cite{DNP2006}, the choice of $h_{d,s}$ has a negligible impact on the the estimators \eqref{f1} and \eqref{f2}   (and the corresponding test)
as long as it  is chosen  sufficiently small.
As a rule of thumb, we choose $h_{d,s}$ as $n_s^{-1/3}$.

For the estimation of the  the long-variance 
we define for $s=1,2$ the partial sum $S_{k,r,s}=\sum_{i=k}^rY_{i,s}$, for some  $m\geq 2$ 
$$
\Delta_{j,s}=\frac{S_{j-m+1,j,s}-S_{j+1,j+m,s}}{m},
$$
and for $t\in[m/n,1-m/n]$
\begin{align}\label{2018-5.4}
\hat \sigma^2_s(t)=\sum_{j=1}^n\frac{m\Delta_{j,s}^2}{2}\omega(t,j),\ \ s=1,2,
\end{align} where for some bandwidth $\tau_{n,s}\in(0,1)$, $$
\omega(t,i)=H\Big(\frac{i/n_s-t}{\tau_{n,s}}\Big)/\sum_{ {i}=1}^nH\Big(\frac{i/n_s-t}{\tau_{n,s}}\Big).
$$
Here $H$ is a symmetric kernel function with compact support $[-1,1]$ and $\int H(x)dx=1$.
For $t\in[0,m/n_s)$ and   $t\in(1-m/n_s,1]$ we define
$\hat \sigma_s^2(t)=\hat \sigma_s^2(m/n_s)$  and 
$\hat \sigma^2(t)=\hat \sigma^2(1-m/n_s)$, respectively.
The consistency of these estimators has been shown in Theorem 4.4 of \cite{dette2018change}.

\subsection{Bootstrap}

Although Theorem \ref{thm32}  is interesting from a theoretical point of view, it cannot be easily implemented for testing the  hypothesis   \eqref{Nullhypo}. The asymptotic bias and variance depend on the long run variances $\sigma_1^2$, $\sigma_2^2$ and the first and second derivative  of the regression functions  $m_1(\cdot)$ and $m_2(\cdot)$. In general, these quantities  are difficult to estimate. Furthermore,
it is well known, that  - even in the case of independence - the convergence rate of statistics as considered  in Theorem \ref{thm32} is slow
(note that  the bias in Theorem \ref{thm32} 
is of order $1/\sqrt{b_{n,1}}$). 
As an alternative
we therefore propose a bootstrap test which does not
require the estimation of the derivatives and addresses the  problem of slow convergence rate. 

The bootstrap procedure is motivated by technical  arguments used in the proof of Theorem \ref{thm32}
in Section \ref{sec6}. There  we show 
(see equations \eqref{2019-32} and \eqref{new.620})
that under the  null hypothesis, the statistic 
$T_{n_1,n_2}$ can be approximated by the  statistic
$$
\int_{\mathbb R} U_n^2(t) w(t)dt,
$$
where 
\begin{align*}
U_n(t)&=\frac{1}{nNb_{n,1}^2h_{d,1}^2}\sum_{j=1}^{n_1}\sum_{i=1}^NK^\circ\Big(\frac{j/n_1-i/N}{b_{n,1}}\Big)K_d'\Big(\frac{ m_1'(i/N)-t}{h_{d,1}}\Big) \sigma_{1}\Big(\frac{j}{n_1}\Big)V_{j,1}\\
&- \frac{1}{nNb_{n,2}^2h_{d,2}^2}\sum_{j=1}^{n_2}\sum_{i=1}^NK^\circ\Big(\frac{j/n_2-i/N}{b_{n,2}}\Big)K_d'\Big(\frac{ m_2'(i/N)-t}{h_{d,2}}\Big) \sigma_{2}\Big(\frac{j}{n_2}\Big)V_{j,2}
\end{align*}
and  $\{V_{j,1},j\in \mathbb Z\}, \{V_{j,2}, j\in \mathbb Z\}$, are sequences of  independent  standard normal distributed random variables.

\begin{algo}
{\rm 
~~

\noindent
(a) Estimate $m_1'$ and $m_2'$ by \eqref{2019-11} and estimate the  long run variances $\sigma^2_1$ and $\sigma^2_2$ by   \eqref{2018-5.4}.
\medskip

\noindent  (b) Generate $B$ copies of standard normal distributed random variables $\{V_{j,1}^{(B)}\}_{j=1}^{n_1}$, $\{V_{j,2}^{(B)}\}_{j=1}^{n_2}$ and calculate the statistic
\begin{align*}
    W_{B}=\int_{\mathbb R} 
    \Big(\frac{1}{nNb_{n,1}^2h_{d,1}^2} \Xi_1^{(B)} (t)-\frac{1}{nNb_{n,2}^2h_{d,2}^2}\Xi_2^{(B)} (t)\Big)^2w(t)dt, 
\end{align*}
where
\begin{align}
    \Xi_1^{(B)}(t)=\sum_{j=1}^{n_1}\sum_{i=1}^NK^\circ\Big(\frac{j/n_1-i/N}{b_{n,1}}\Big)K_d'\Big(\frac{\hat m_1'(i/N)-t}{h_{d,1}}\Big)\hat \sigma_{1}\Big(\frac{j}{n_1}\Big)V_{j,1}^{(B)},\\
    \Xi_2^{(B)}(t)=\sum_{j=1}^{n_2}\sum_{i=1}^NK^\circ\Big(\frac{j/n_2-i/N}{b_{n,2}}\Big)K_d'\Big(\frac{\hat m_2'(i/N)-t}{h_{d,2}}\Big)\hat \sigma_{2}\Big(\frac{j}{n_2}\Big)V_{j,2}^{(B)}.
\end{align}
\medskip
\noindent 
 (c) Let $W_{(1)}\leq W_{(2)}\leq\ldots \leq W_{(B)}$ be the ordered statistics of $\{W_{s},1\leq s\leq B\}$. We reject 
  the null hypothesis \eqref{Nullhypo}
  at level $\alpha$, whenever
  \begin{equation}
      \label{boottest}
       T_{n_1,n_2} > W_{(\lf B(1-\alpha)\rf)}. 
  \end{equation}
  The $p$-value  of this  test is given by $1-B^*/B$.
where $B^*=\max\{r:W_{(r)}\leq T_{n_1,n_2}\}$. 
}
\end{algo}

\subsection{Simulated level and power}
\label{Sec4.1}

In this section we illustrate the finite sample properties of the test \eqref{boottest} by means of a small simulation study. All presented results
are based on   $1000$ runs
and $B=500$ bootstrap replications. 
We consider equal  sample sizes  $n_1 = n_2 = n = 100$, $200$ and $500$. Throughout this article, the  Epanechnikov kernel (e.g., see \cite{Silverman1998}) is considered for all kernels
appearing in the test procedure, and  we use $N = n$ in \eqref{f1} and \eqref{f2}.
Besides, $h_{d, N} = n^{-1/3}$, and $b_{n_1}$ and $b_{n_2}$ are  chosen as described at  the beginning of Section \ref{sec4}.

For $s = 1$ and $2$, we consider model \eqref{Model_1} with the error process
\begin{align}
\label{error}
G_s(t,\FF_i)=0.6(t-0.3)^2G(t,\FF_{i-1,s})+\eta_{i,s},
\end{align}
where $\FF_{i,s}=(...,\eta_{i-1,s},\eta_{i,s})$. We assume that $\eta_{i,1}$ are i.i.d standard normal random variables, and $\eta_{i,2}$ are i.i.d. copies of the random variable $t_{5}/\sqrt{5/3}$, where $t_{5}$ denotes the  $t$-distribution with $5$ degrees of freedom.   For the regression functions we consider the models
\begin{align}
\label{exA} m_{1} (x) &=  (x - 0.4)^{2} \text{ ~~and  ~} m_{2} (x) = (x - 0.3)^{2} - 0.2 ,  \\
\label{exB}
 m_{1}(x) &=  \sin (-\pi x)  \text{ ~~and  ~ } m_{2}(x) =\sin(-\pi(x + 0.1)) + \frac{1}{4}.
\end{align}
In  Table \ref{tab3} we display the rejection probabilities of the test \eqref{boottest}, where the level of significance is $5\%$ and $10\%$. 
The results show a good approximation of the nominal level
in all cases under consideration.  

\begin{table}[h!]

\begin{center}

\begin{tabular}{|c|c|c|c|}\hline
model  & $n = 100$ & $n = 200$ & $n = 500$\\ \hline 
\eqref{exA} & {$0.057$} & {$0.054$} & {$0.051$}\\ \hline
\eqref{exB} & {$0.059$} & {$0.057$} & {$0.054$}\\ \hline
\hline
\eqref{exA} & { 0.111} & { 0.108} & {0.104} \\ \hline
\eqref{exB} & { 0.116} & { 0.112} & { 0.103} \\ \hline
\end{tabular}
\end{center}
\caption{\it   The estimated size of the test  \eqref{boottest} for different sample sizes $n_1= n_2=n$. The level of significance is  $5\%$ (upper part) and $10\%$ (lower part).}

\label{tab3}

\end{table}

 In order to study the power of the test 
 \eqref{boottest}  we consider 
 the same error processes as in \eqref{error} 
 and used  the regression functions

    \begin{align}\label{exC} 
 m_{1} (x) & =  (x - 0.4)^{2} \text{~~ and~ } m_{2} (x) = x^{3}, \\ 
 m_{1}(x) &= \sin(-\pi x) \text{~~ and~ } 
 m_{2}(x) = -\cos(\pi x). 
 \label{exD} 
\end{align}
The simulated power  is displayed in Table \ref{tab4} and the  results indicate that the test 
detects  the alternatives reasonably well.

\begin{table}[h!]

\begin{center}

\begin{tabular}{|c|c|c|c|}\hline
model  & $n = 100$ & $n = 200$ & $n = 500$\\ \hline 
\eqref{exC}  & { $0.563$} & { $0.647$} & { $0.778$}\\ \hline
\eqref{exD}  & { $0.617$} & { $0.744$} & { $0.822$}\\ \hline
 \hline  \hline 
\eqref{exC} & { $0.722$} & {0.847} & {0.899} \\ \hline
\eqref{exD} & { $0.777$} & { 0.868} & {0.971} \\ \hline
\end{tabular}
\end{center}

\caption{\it   The estimated power  of the test  \eqref{boottest} for different sample sizes $n_1= n_2=n$. The level of significance is  $5\%$ (upper part) and $10\%$ (lower part).}

\label{tab4}

\end{table}

\begin{figure}[p]
\includegraphics[width= 3.5in, height = 3in]{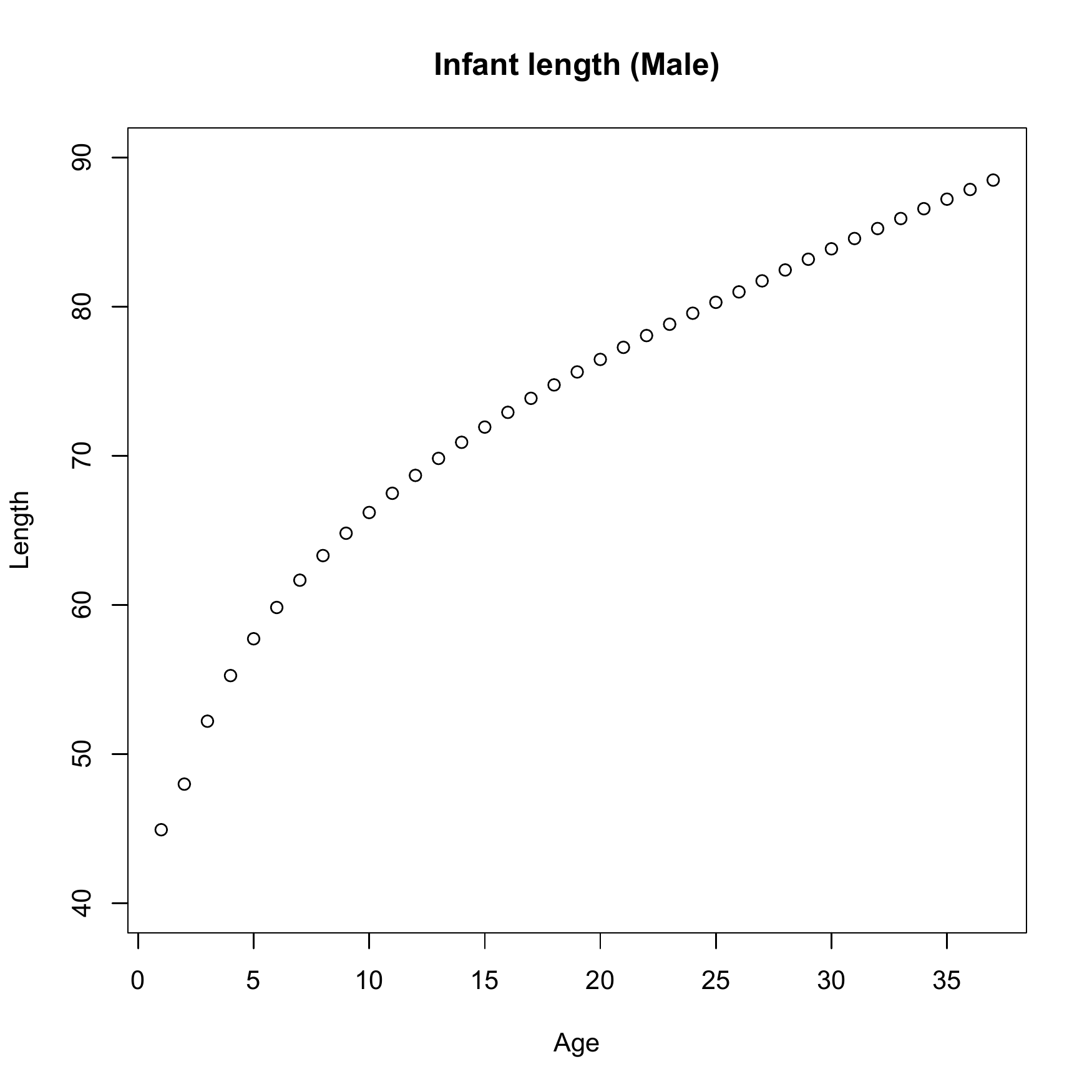}\includegraphics[width= 3.5in, height = 3in]{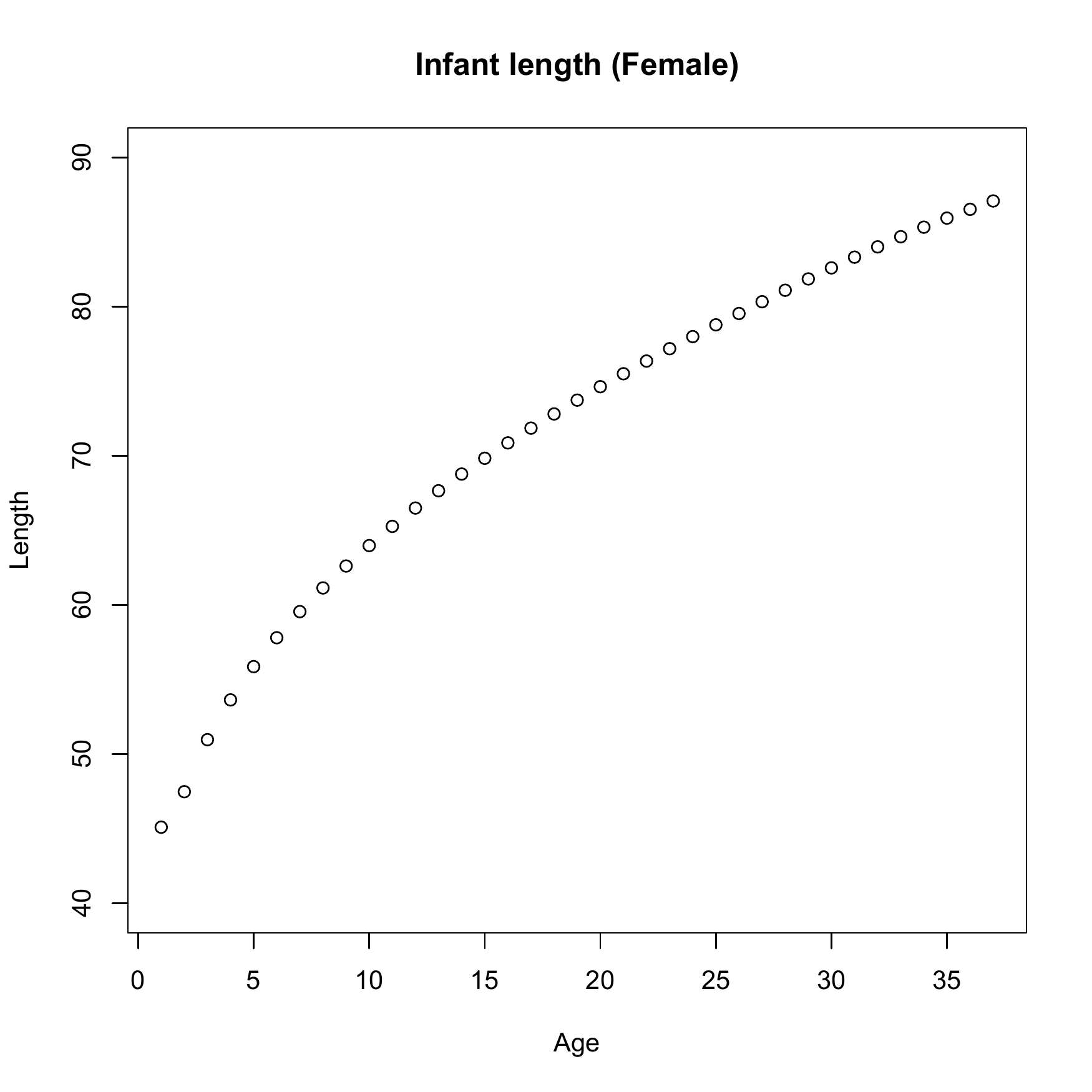}\\ 
\caption{\it  
\label{Figreal}
Plots of the length of the male (right part) and female (left parts) infants for different age.
 } 
\end{figure}

\begin{figure}
\begin{center}
\includegraphics[width= 3.5in, height = 3in]{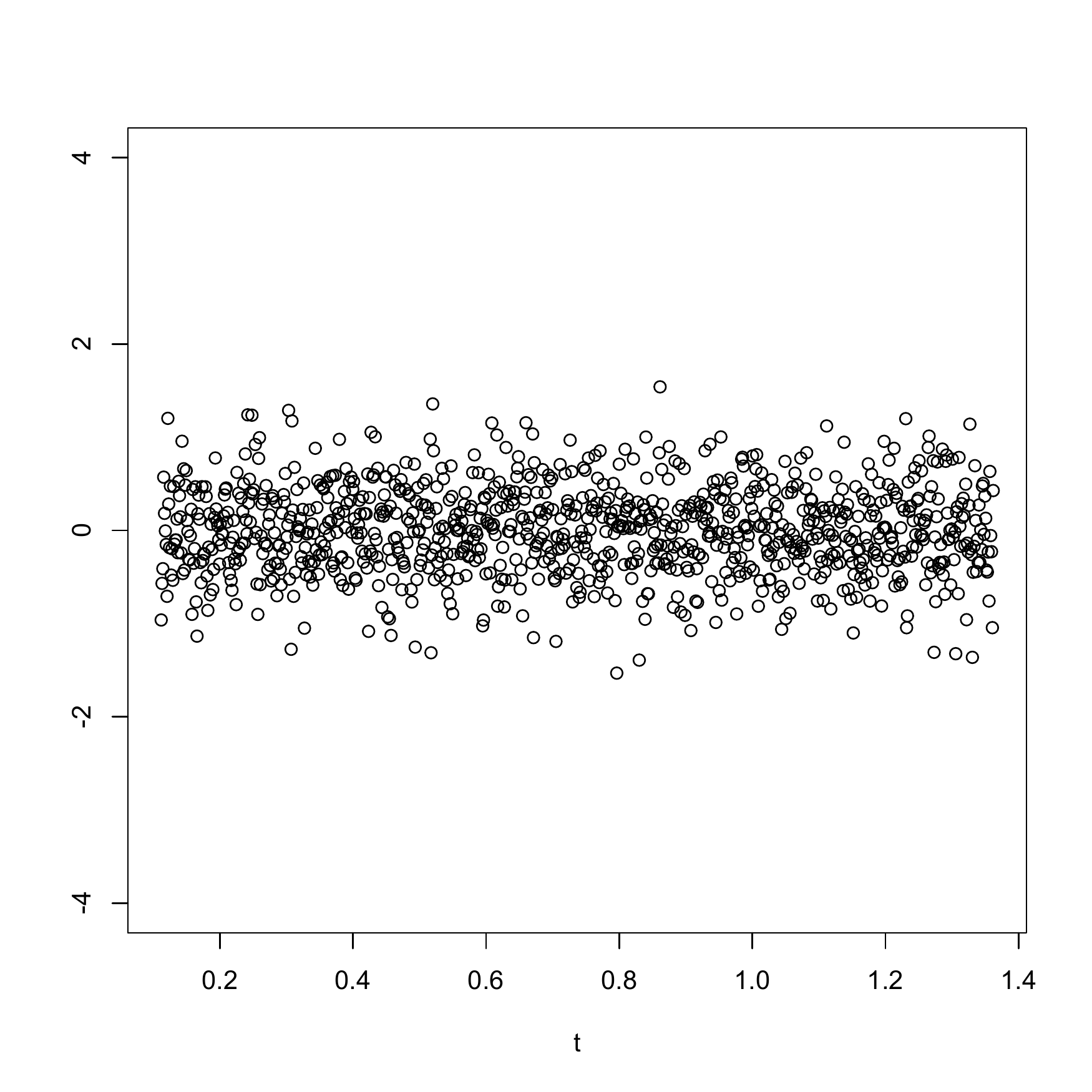}\\ 
\caption{\it  
\label{GraReal}
Plots of ${\cal{C}}_{n_{1}, n_{2}}$ for the real data described in Section 4.3.} 
\end{center}
\end{figure}

\subsection {Real data analysis}
In this section, we use   the test \eqref{boottest} and the 
graphical device described in Section \ref{sec21} to 
investigate the validity of assertion \eqref{Nullhypo}  for  growth data of  male and female infants. This data set is available from  \url{https://www.cdc.gov/growthcharts/html_charts/lenageinf.htm#males} and consists  of the  monthly growth of length of 
male and female infants in the first three years (here $n_{1} = n_{2} = 37$). The data is depicted in Figure \ref{Figreal} and indicates that the relation between length and age in both groups might be concave. Therefore  we
model the negative values  of this data  by  two regression  models  of the form  \eqref{Model_1} with convex  regression functions, where   group 1 represents the male  and group 2 the female infants.  For this data, we obtain $\hat{c} = 0.046$
as estimate for the  horizontal shift using the 
statistic  \eqref{chat} 
and $\hat{d} = \hat{m}_{1}(0) - \hat{m}_{2}(\hat{c})= 0.087$ 
as  estimate  of the vertical  shift $d$. 

We begin illustrating the application of the graphical device 
described in Section \ref{sec21}. In Figure \ref{GraReal} we plot the points of  the set  ${\cal{C}}_{n_1, n_{2}}$  in \eqref{cn}  using   $L = 1000$ equally spaced points in the interval  $(\hat{a} + \eta, \hat{b} - \eta)$, where $\hat{a} = \hat{m}_{1}^{'} (0) = 0.112$, $\hat{b} = \hat{m}_{1}^{'} (1 - \hat{c}) = 1.362$, and $\eta = 0.001$ is chosen (the smoothing parameters are chosen as described in Section \ref{sec4}). The figure  clearly indicates the existence of a vertical and horizontal  shift between the regression functions as formulated in the  null hypothesis  \eqref{Nullhypo}.

Finally, we also investigate the performance of the  test 
\eqref{boottest} for this  data set, where  all parameters required  for the bootstrap test  are chosen as described in Section \ref{sec4}. For $B = 500$ bootstrap replications, we obtain the $p$-value  $0.781$, which gives no indication to reject
the null hypothesis and is   consistent with the conclusion made by graphical inspection.

\label{sec5}
\def\theequation{5.\arabic{equation}}
\setcounter{equation}{0}

\section{Appendix : Proofs} 
\def\theequation{5.\arabic{equation}}
\setcounter{equation}{0}


\label{sec6}
\subsection{Preliminaries} \label{sec61}

In this section, we state a few auxiliary results, which will be used later in the proof.
We begin with Gaussian approximation. A proof of this result can be found in \HE{\cite{wu2011gaussian}}.

\begin{proposition}\label{Prop3.1} 
Let
$$
\bs S_i=\sum_{s=1}^i \mathbf e_i,
$$
and assume that the Assumption \ref{Assum3.1} is satisfied.
	Then on a possibly richer probability space, there exists a process  $\{\mathbf S_i^\dag\}_{i\in \mathbb Z}$ such that
	$$
	\{\bs S_i^\dag\}_{i=0}^n
	\stackrel{\cal D}=\{\bs S_i\}_{i=0}^n
	$$
	(equality in distribution), and a sequence of independent $2$-dimensional standard normal 
distributed 	random variables $\{\mathbf V_i\}_{i\in \mathbb Z}$, such that
	\begin{align}
	\max_{1\leq j\leq n}\Big |\sum_{i=1}^j\mathbf S_i^\dag-\sum_{i=1}^j \Sigma(i/n)\mathbf V_i\Big |= o_p(n^{1/4}\log^{2} n),
	\end{align}
	where $\Sigma(t)$ is the square root of the long-run variance  matrix $\Sigma^2(t)$
defined in Assumption \ref{Assum3.1}.
	 \end{proposition}
	
\begin{proposition}\label{PropM}
	Let Assumption \ref{Assum3.1} and \ref{AssumptionK} be satisfied.
\begin{itemize}
    \item[(i)]   For $s=1,2$ we have
  \begin{align}\label{Jul1-9}
  \sup_{ t\in [b_{n,s},1-b_{n,s}]}\Big | \hat m_s'(t)-m_s'(t)-\frac{1}{n_sb_{n,s}^2}\sum_{i=1}^{n_s}K^\circ\Big(\frac{i/n_s-t}{b_{n,s}}\Big)e_{i,s}\Big |=O_P\Big(\frac{1}{n_sb_{n,s}^2}+b_{n,s}^2\Big)
\end{align}
 where the kernel $K^\circ$ is defined in \eqref{kequiv}.  
 \item[(ii)]
 For  $s=1,2$ 
	\begin{align}\label{Jul1-10}
	\sup_{t\in[b_{n,s},1-b_{n,s}]}\Big|\frac{1}{n_sb_{n_s}^2}\sum_{i=1}^{n_s} K^\circ\Big(\frac{i/n_s-t}{b_{n,s}}\Big)\Big(e_{i,s}-\sigma_s(i/n) V_{i,s}\Big)\Big|=o_p\Big(\frac{\log ^2n_s}{n^{3/4}_sb_{n,s}^2}\Big),
	\end{align} 
	where $\{ V_{i,s}, ~ i=1, \ldots , n_s ,s=1,2\}$ 
	denotes a sequence of independent  standard normal distributed random random variables.
	\item[(iii)] 
	For $s=1,2$ we have
	\begin{align}\label{supM}
	\sup_{t\in [b_{n,s},1-b_{n,s}]}|\hat m'_s(t)-m_s'(t)|=O_p\Big(\frac{\log n_s}{\sqrt{n_sb_{n,s}}b_{n,s}}+\frac{\log ^2n_s}{n^{3/4}_sb_{n,s}^2}+b_{n,s}^2\Big).
		\end{align} 
	
		\item[(iv)] 
	For $s=1,2$ we have
	\begin{align}\label{supM0}
	\sup_{t\in [0,b_{n,s}]\cup[1-b_{n,s},1] }|\hat m'_s(t)-m_s'(t)|=O_p\Big(\frac{\log n_s}{\sqrt{nb_{n,s}}b_{n,s}}+\frac{\log ^2n_s}{n_s^{3/4}b_{n,s}^2}+b_{n,s}\Big).
		\end{align} 
	\end{itemize}
\end{proposition}

\medskip

\noindent {\bf Proof:} 

\noindent 
(i):
Define for $s=1,2$ and $l=0,1,2$
\begin{align*}
    R_{n,s,l}(t) & =\frac{1}{n_sb_{n,s}}\sum_{i=1}^{n_s}Y_{i,s}K\Big(\frac{i/n-t}{b_{n,s}}\Big)\Big(\frac{i/n_s-t}{b_{n,s}}\Big)^l, \\
     S_{n,s,l}(t) & =\frac{1}{n_sb_{n,s}}\sum_{i=1}^{n_s}K\Big(\frac{i/n_s-t}{b_{n,s}}\Big)\Big(\frac{i/n_s-t}{b_{n,s}}\Big)^l
\end{align*}
 Straightforward calculations show that 
\begin{align*}
    (\hat m_s(t), b_{n,s}\hat m_s'(t))^\top=S_{n,s}^{-1}(t)R_{n,s}(t)
   ~~~~~~(s=1,2) ,    \end{align*}
 where
 \begin{align*}
  R_{n,s}(t)=
  \begin{pmatrix} 
  R_{n,s,0}(t) \\ 
  R_{n,s,1}(t)
  \end{pmatrix}
  ~  ,~~ S_{n,s}(t)=\begin{pmatrix} 
S_{n,s,0} & S_{n,s,1}  \\
S_{n,s,1} & S_{n,s,2}.
\end{pmatrix}
   \end{align*}
Note that   Assumption \ref{AssumptionK} 
gives 
\begin{align}
    S_{n,s,0}(t)=1+O\Big(\frac{1}{n_sb_s}\Big),S_{n,s,1}(t)=O\Big(\frac{1}{n_sb_{n,s}}\Big),S_{n,s,2}(t)=\int_{-1}^1K(x)x^2dx+O\Big(\frac{1}{n_sb_{n,s}}\Big)
\end{align}
 uniformly with respect to $t\in [b_{n,s},1-b_{n,s}]$. 
 The first part of the proposition  now follows by a  Taylor expansion of $R_{n,s,l}(t)$.

\medskip 
\noindent
(ii):  The fact asserted in \eqref{Jul1-10} follows from \eqref{Jul1-9}, Proposition \ref{Prop3.1}, the  summation by parts formula and similar arguments to derive equation (44) in \cite{zhou2010nonparametric}. 

\medskip\noindent
(iii) + (iv):  Following Lemma 10.3 of \cite{dette2018change}, we have
\begin{align}\label{Jul1-11}
\sup_{t\in [b_{n,s},1-b_{n,s}]}\Big|\frac{1}{n_sb_{n,s}}\sum_{i=1}^{n_s}K^\circ\Big(\frac{i/n_s-t}{n_sb_{n,s}}\Big)\Big(\sigma_s(\frac{i}{n_s})V_{i,s}\Big)\Big|=O_p\Big(\frac{\log n_s}{\sqrt{n_sb_{n,s}}}\Big).
\end{align}

Finally, \eqref{supM} and \eqref{supM0} follow from \eqref{Jul1-9} \eqref{Jul1-10} and  \eqref{Jul1-11}, which completes the proof of Proposition \ref{PropM}. 
\hfill $\Box$

\subsection{Proof of Theorem \ref{The0}} 
\label{sec61a}
We only prove the result in the case $g \equiv 0 $. The general case follows by  the same  arguments.
Under Assumptions \ref{Assum3.1} and \ref{AssumptionK}, it follows from the proof of Theorem 4.1 in \cite{dette2018change}  that 
$$
\displaystyle\sup_{t\in(a+\eta, b - \eta)} \big[\big(\hat f_1(t) - \hat f_2(t)\big) - \big(((m_1')^{-1}(t))' - ((m_{2}')^{-1})'(t)\big)\big]\rightarrow 0
$$ in probability, where 
$\hat f_1^{-1} (t)$ and $\hat f_{2}(t)$ are  defined in \eqref{f1} and \eqref{f2}, respectively.
 Next, since under the null hypothesis \eqref{Nullhypo}, $((m_1')^{-1}(t))' - ((m_{2}')^{-1})'(t) = 0$ for all $t\in (a+\eta, b - \eta)$, (See Lemma \ref{lem1}) we have under the null hypothesis, $$\displaystyle\sup_{t\in(a+\eta, b - \eta)} \big [\hat f_1(t) - \hat f_2(t)  \big ]\rightarrow 0$$ in probability. In other words, under $H_{0}$, for any $\epsilon > 0$, we have \begin{align}
    \displaystyle\lim_{n\rightarrow\infty}
    \mathbb P 
    \Big [\displaystyle\sup_{t\in(a+\eta, b - \eta)}\big|\hat f_1(t)- \hat f_2(t)\big | <\epsilon\Big ] = 1, 
\end{align} and hence, under the null hypothesis  $g \equiv 0$,  we have  $\mathbb P [{\cal{C}}_{n_1,n_2}\subset L(\epsilon)] = 1$. 
\hfill$\Box$

\noindent 
\subsection{Proof of Theorem \ref{thm32}}
\label{sec62}

To simplify the notation, we prove Theorem \ref{thm32} in  the case of equal sample sizes
and equal bandwidths.  The general case follows by  the same arguments with an additional amount of notation. In this case $c_2=r_2=1$ and we omit the subscript in bandwidths if no confusion arises, for example we write  $n_1=n_2=n$,   $b_{n,1}=b_{n,2} =b_n$
and use a similar notation for other symbols depending on the sample size. In particular, we write $T_n$ for $T_{n_1,n_2}$ if $ n= n_1=n_2$.

Define the statistic 
\begin{align}
    \tilde T_n=\int\Big(\hat f_1(t)-\hat f_2(t)\Big)^2 w(t)dt
\end{align}
which is obtained from $T_n$ by replacing the weight function $\hat w$ in \eqref{def_T} by its deterministic analogue \eqref{aandb}.
We shall show Theorem \ref{thm32} in two steps proving the assertions
\begin{align}
\label{s1}
&& nb_n^{9/2}\tilde T_n- B_n(g) \Rightarrow {\cal N} (0, V_T) \\
&&\label{s2}
nb_n^{9/2}(T_n-\tilde T_n)=o_p(1).
\end{align}

\subsubsection{Proof of  \eqref{s1}} 
By  simple algebra, we obtain the decomposition
\begin{align}
\tilde T_n=\int (I_1(t)-I_2(t)+II(t))^2w(t)dt, 
\end{align} 
where for $s=1,2$
\begin{align}
I_s(t)=\frac{1}{Nh_d}\sum_{i=1}^N\Big(K_d\Big(\frac{\hat m_s'(i/N)-t}{h_d}\Big)-K_d\Big(\frac{ m_s'(i/N)-t}{h_d}\Big)\Big),
\label{is}\\
II(t)=\frac{1}{Nh_d}\sum_{i=1}^N\Big(K_d\Big(\frac{ m_1'(i/N)-t}{h_d}\Big)-K_d\Big(\frac{ m_2'(i/N)-t}{h_d}\Big)\Big).
\label{ii}
\end{align}
Observing the estimate on page  471  of   \cite{DNP2006} it follows
\begin{align}
\frac{1}{Nh_d}\sum_{i=1}^NK_d\Big(\frac{m_s'(i/N)-t}{h_d}\Big)=\Big(((m_s')^{-1}(t))'+O\Big(h_d+\frac{1}{Nh_d}\Big)\Big)~,
\end{align}
($s=1,2$) which yields  the estimate
\begin{align}\label{Jul-21}
II(t)=((m_1')^{-1}(t))'-((m'_2)^{-1}(t))'+O\Big(h_d+\frac{1}{Nh_d}\Big)
\end{align}
uniformly with respect to  $t\in [a+\eta,b-\eta]$.
 For  the two other  terms we use a Taylor expansion
 and obtain the decomposition 
\begin{align} \label{tdec} 
I_s(t)=I_{s,1}(t)+I_{s,2}(t)  ~~~~(s=1,2) ,
\end{align}
where 
\begin{align*}
&I_{s,1}(t)=\frac{1}{Nh^2_d}\sum_{i=1}^N
K'_d\Big(\frac{m_s'(i/N)-t}{h_d}\Big)(\hat m_s'(i/N)-m_s'(i/N)),\\
&
I_{s,2}(t)=\frac{1}{2Nh^3_d}\sum_{i=1}^{N}
K''_d\Big(\frac{m_s'(i/N)-t+\theta_s(\hat m_s'(i/N)-m_s(i/N))}{h_d}\Big)(\hat m_s'(i/N)-m_s'(i/N))^2
\end{align*}
for some $\theta_s\in [-1,1]$ ($s=1,2$).
By part (iii) and (iv) of Proposition \ref{PropM} and the  same  arguments that were used in the online  supplement of  \cite{dette2018change},
to obtain the bound for the  term $\Delta_{2,N}$ in the proof of their Theorem 4.1 
 it follows that 
\begin{align}\label{Jul-24}
I_{s,2} (t)=O_p\Big (\frac{\pi_n^2}{h_d^3}(h_d+\pi_n)\Big)=O_p\Big(\frac{\pi_n^3}{h_d^3}\Big) ~(s=1,2), 
\end{align}
uniformly with respect to $t\in [a+\eta,b-\eta]$.
Here we used the fact that  the number of non-zero summands in $I_{s,2}(t)$ is  of order $O(h_d+\pi_n)$.

Next, for the investigation  of the difference $I_{1,1}(t)-I_{2,1}(t)$,  we define $\mathbf m'= (m_1,m_2) $ and consider the vector
$$K'_d\Big( \frac{\mathbf m'(i/N)-t}{h_d}\Big)=\Big (K'_d\Big( \frac{ m'_1(i/N)-t}{h_d}\Big),-K'_d\Big( \frac{ m'_2(i/N)-t}{h_d}\Big)\Big )^\top.$$
By part  (i)  and (ii) of  Proposition
\ref{PropM},  it follows that there exists 
independent $2$-dimensional  standard normal distributed 
random vectors   $\mathbf V_i $  such that
\begin{align*}
I_{1,1} (t)-I_{2,1}(t) & =\frac{1}{nNb_n^2h^2_d}\sum_{j=1}^n\sum_{i=1}^NK^\circ\Big(\frac{j/n-i/N}{b_n}\Big)(K'_d)^{T}\Big(\frac{\mathbf m'(i/N)-t}{h_d}\Big)\Sigma(j/n)\mathbf V_j \\
&~~~~~~~~~~~~~~~~~~~~~~~~~~~~~
~~~~~~~~
+{ O_p(\pi_n'h_d^{-1})}.
\end{align*}
uniformly with respect to  $t\in [a+\eta,b-\eta]$.
Combining  this estimate  with equations \eqref{Jul-21} and \eqref{Jul-24}, it follows
\begin{align}\label{2019-32}
T_n=\int \left(U_n(t)+((m'_1)^{-1}(t))'-((m'_2)^{-1}(t))'+R_n^\dag(t)\right)^2w(t)dt, 
\end{align}
where
\begin{align}
U_n(t) & =\frac{1}{nNb_n^2h^2_d}\sum_{j=1}^n\sum_{i=1}^NK^\circ\Big(\frac{j/n-i/N}{b_n}\Big)(K'_d)^\top\Big(\frac{\mathbf m'(i/N)-t}{h_d}\Big)\Sigma(j/n)\mathbf V_j,\label{new.620}
\end{align}
and the  remainder $R_n^\dag(t)$ can be estimated as follows
\begin{align}
\sup_{t\in [a+\eta,b-\eta]}|R^\dag_n(t)|
& ={O_p\Big(\frac{\pi_n'}{h_d}+\frac{\pi_n^3}{h_d^3}+h_d+\frac{1}{Nh_d}\Big)}\label{new.28}.
\end{align}
 We now study the asymptotic properties of 
 to the quantities 
\begin{align}\label{Un(t)}
& nb_n^{9/2}\int (U_n(t))^2w(t)dt, \\
& \label{Un(t)1}
nb_n^{9/2}\int U_n(t)((m_1^{-1}(t))'-(m_2^{-1}(t))')w(t)dt, \\
& nb_n^{9/2}\int U_n(t)R_n^\dag(t)w(t)dt,
\label{Un(t)2}
\end{align} 
which  determine the asymptotic distribution of $T_{n}$ since the bandwidth conditions yield under local alternatives 
in the case  $(m_1^{-1}(t))'-(m_2^{-1}(t))'=\rho_ng(t)$,  
\begin{align} \label{a1}
   nb_n^{9/2}\int \rho_n^2(t)w(t)=\int g^2(t)w(t)dt,
\end{align}
and the other parts of the expansion are negligible, i.e., 
\begin{align}
\label{a2}
& nb_n^{9/2}\int (R_n^\dag(t))^2w(t)dt=o(1), \\ 
&  nb_n^{9/2}\int \rho_ng(t)R^\dag_{n}(t)w(t)dt=o(1).
\label{a3}
\end{align}

\paragraph{Asymptotic properties  of \eqref{Un(t)}:}
To address the expressions related to $U_n(t)$ in \eqref{Un(t)} -  \eqref{Un(t)2} note that 
\begin{align}
U_n(t)=U_{n,1}(t)-U_{n,2}(t),
\end{align}
where 
\begin{align}
U_{n,s}(t)=\frac{1}{nNb_n^2h_d^2}\sum_{j=1}^n\sum_{i=1}^NK^\circ\Big(\frac{j/n-i/N}{b_n}\Big)K'_d\Big({ m_s'(i/N)-t}{h_d}\Big)\sigma_s(j/n) V_{j,s}
\end{align}
for $s=1,2$, and $\{V_{j,s}\}$ are  independent  standard normal distributed random variables. 
In order to simplify the notation, we define the quantities
$$
U_{n,s}(t)=\displaystyle\sum_{j=1}^nG(m_s'(\cdot),j,t)V_{j,s}~~  ~~  ~~  ~~  ~~  (s=1,2),
$$
where
\begin{align}
G(m_s'(\cdot),j,t)=\frac{1}{nNb_n^2h_d^2}\sum_{i=1}^NK^\circ\Big(\frac{j/n-i/N}{b_n}\Big)K'_d\Big(\frac{ m_s'(i/N)-t}{h_d}\Big)\sigma_s(j/n).
\end{align}
A straightforward calculation 
(using the  change of variable
$v=(m_s'(u)-t)/h_d$)
shows that  
\begin{align*}
G( m'_s(\cdot),j,t)&=\frac{1}{nb_n^2h^2_d}\int_{0}^{1}K^\circ\left(\frac{j/n-u}{b_n}\right)K'_d\left(\frac{ m_s'(u)-t}{h_d}\right)\sigma_s(j/n) du+O\left(\delta_n 
\right) 
\notag\\
\notag
&=\frac{1}{nb_n^2h_d}\sigma_s(j/n)\int_{\mathcal A_s(t)}K'_d(v)((m_s')^{-1}(t+h_dv))'K^\circ\left(\frac{j/n-(m_s')^{-1}(t+h_dv)}{b_n}\right)dv\notag\\
&+O\left(\delta_n 
\right), 
\end{align*}
where the interval $\mathcal A_s(t)$ is defined by 
\begin{align}
\mathcal A_s(t)=\Big(\frac{m_s'(0)-t}{h_d},\frac{ m_s'(1)-t}{h_d}\Big),
\end{align}
the remainder is given by 
$$
\delta_n = O\Big(\Big(\frac{1}{nb_n^2h^2_dN}\Big)\mathbf 1\Big(\Big|\frac{j/n-(m_s')^{-1}(t)}{b_n+Mh_d}\Big|\leq 1\Big)\Big), 
$$
and $\mathbf 1 (A) $  denote the indicator function of the set $A$.
As the kernel $K'_d(\cdot)$ has a compact support and is symmetric, it  follows by  a Taylor expansion
 for  any  $t$  with  $w(t)\neq 0$
\begin{align*}
& \int_{\mathcal A_s(t)}K'_d(v)((m_s')^{-1}(t+h_dv))'K^\circ\Big(\frac{j/n-(m_s')^{-1}(t+h_dv)}{b_n}\Big)dv\notag\\   
& =-\frac{h_d}{b_n}(((m_s')^{-1}(t))')^2(K^\circ)'\Big(\frac{j/n-(m_s')^{-1}(t)}{b_n}\Big)\int K_d'(v)vdv\Big(1+O\Big(b_n+\frac{h_d^2}{b_n^2}\Big)\Big)
\end{align*} 
With the notation 
\begin{align}
\tilde G(m_s'(\cdot),j,t)=\frac{-1}{nb_n^3}(K^\circ)'\Big(\frac{j/n-(m_s')^{-1}(t)}{b_n}\Big)\sigma_s(j/n)(((m_s')^{-1})'(t))^2\int vK'_d(v)dv
\end{align}
($s=1,2$) we thus obtain the approximation 
\begin{align}\label{decom1}
\int U_n^2(t)w(t)
&=\sum_{s=1}^2 \sum_{j=1}^nV_{j,s}^2\int G^2(m_s'(\cdot),j,t)^2w(t)dt
\\&+\sum_{s=1}^2\sum_{1\leq i\neq j \leq n}V_{i,s}V_{j,s}\int  G(m_s'(\cdot),i,t) G(m_s'(\cdot),j,t)w(t)dt\notag
\\&-2\sum_{1\leq i\leq n}V_{i,1}V_{i,2}\int G(m_1'(\cdot),i,t) G(m_2'(\cdot),i,t)w(t)dt\notag
\\&=\sum_{s=1}^2 \sum_{j=1}^nV_{j,s}^2\left(\int\tilde G^2(m_s'(\cdot),j,t)^2w(t)dt(1+r_{i,s})\right)\notag
\\&+\sum_{s=1}^2\sum_{1\leq i\neq j \leq n}V_{i,s}V_{j,s}\left(\int \tilde  G(m_s'(\cdot),i,t)\tilde G(m_s'(\cdot),j,t)w(t)dt(1+r_{i,j,s})\right)\notag
\\&-2\sum_{1\leq i\leq n}V_{i,1}V_{i,2}\left(\int \tilde G(m_1'(\cdot),i,t)\tilde G(m_2'(\cdot),i,t)w(t)dt(1+r_{i,s}')\right)\notag,
\end{align}
where the remainder satisfy 
$$
\max\left(\displaystyle\max_{i,j,s=1,2}(|r_{i,j,s}|),\displaystyle\max_{i,s=1,2}(|r_{i,s}|),\max_{i,s=1,2}(|r_{i,s}'|)\right)=o(1).
$$
Let us now consider  the statistcis  $\tilde U_{n,s}(t)=\sum_{j=1}^n\tilde G(m_s'(\cdot),j,t)V_{j,s}$ ($s=1,2$), and \begin{align}\label{tildeUn(t)}
    \tilde U_n(t)=\tilde U_{n,1}(t)-\tilde U_{n,2}(t), 
\end{align} 
then, by the previous calculations, it follows that
\begin{align}\label{equiv}
  nb_n^{9/2} \Big ( \int U_n^2(t) w(t)dt
- \int  \tilde U_n^2(t)w(t)dt \Big ) = o_P(1),
\end{align} 
and therefore, we investigate  the weak convergence 
of the statistic $nb_n^{9/2}  \int  \tilde U_n^2(t)w(t)dt$ in the following.
For this purpose we use a similar decomposition 
as in \eqref{decom1} and obtain 
\begin{align}
\int \tilde U_n^2(t) w(t)dt&=\sum_{s=1}^2\int (\tilde U_{n,s}(t))^2w(t)dt-2\int (\tilde U_{n,1}(t)\tilde U_{n,2}(t))w(t)dt\notag
\\&=\sum_{s=1}^2 \sum_{j=1}^nV_{j,s}^2\int\tilde G^2(m_s'(\cdot),j,t)^2w(t)dt\notag
\\&+\sum_{s=1}^2\sum_{1\leq i\neq j \leq n}V_{i,s}V_{j,s}\int \tilde G(m_s'(\cdot),i,t)\tilde G(m_s'(\cdot),j,t)w(t)dt\notag
\\&-2\sum_{1\leq i\leq n}V_{i,1}V_{i,2}\int \tilde G(m_1'(\cdot),i,t)\tilde G(m_2'(\cdot),i,t)w(t)dt\notag
\\&:=D_1+D_2+D_3,
\label{decom2}
\end{align}
where the last equation defines $D_1,D_2$ and $D_3$ in an obvious manner. 
Elementary calculations (using a Taylor expansion and the fact that the kernels have compact support)
show that
\begin{align}
\E (D_1&)=\sum_{s=1}^2\sum_{j=1}^n\int\Big(\frac{-1}{nb_n^3}(K^\circ)'\Big(\frac{j/n-(m_s')^{-1}(t)}{b_n}\Big)\sigma_s(j/n)(((m_s')^{-1})'(t))^2\int vK'_d(v)dv\Big)^2w(t)dt\notag
\\&=\sum_{s=1}^2\sum_{j=1}^n\int\Big(\frac{1}{nb_n^3}(K^\circ)'\Big(\frac{j/n-(m_s')^{-1}(t)}{b_n}\Big)\sigma_s((m_s')^{-1}(t))(((m_s')^{-1})'(t))^2\int vK'_d(v)dv\Big)^2\notag\\&\times w(t)dt(1+O(b_n))\label{7-21-36}.
\end{align}
 Using  the estimate 
\begin{align}
\frac{1}{nb_n}\sum_{j=1}^n\Big((K^\circ)'\Big(\frac{j/n-(m_s')^{-1}(t)}{b_n}\Big)\Big)^2=\int ((K^\circ)'(x))^2dx\Big(1+O\Big(\frac{1}{nb_n}\Big)\Big),\label{7-21-37}
\end{align}
(uniformly with respect to  $t\in [a+\eta, b-\eta]$), \eqref{7-21-36} and \eqref{7-21-37} gives
\begin{align*}
\E (D_1) &=\frac{1}{nb_n^5}\sum_{s=1}^2\int ((K^\circ)'(x))^2dx\int \left(\sigma_s((m_s')^{-1}(t))(((m_s')^{-1}(t))')^2\int vK'_d(v)dv\right)^2w(t)dt \\
& ~~~~~~~~~~~~~~~~~~~~~~~~
\times
\Big (1+O\Big (b_n+\frac{1}{nb_n}\Big )\Big),
\end{align*}
which  implies 
\begin{align}\label{eq37}
	\E(nb_n^{9/2}D_1)=
	B_n(0) +O\Big(\sqrt{b_n}+\frac{1}{nb_n^{3/2}}\Big),
	\end{align}
	where $B_n (g)$  is defined in Theorem \ref{thm32} (and we use the notation with the function $g \equiv 0$).
Here we  used the  change of   variable $(m_s')^{-1}(t)=u$, and afterwards, 
	$((m_s')^{-1})'(t)=\frac{1}{m_s''((m_s')^{-1}(t))}$.
Similar arguments establish that 
\begin{align*}
\mbox{Var} (D_1)=O\Big(\sum_{s=1}^2 \sum_{j=1}^n(\int\tilde G^2(m_s'(\cdot),j,t)^2w(t)dt)^2\Big)=O\Big(\frac{nb_n^2}{n^4b_n^{12}}\Big)=O\Big(\frac{1}{n^3b_n^{10}}\Big),
\end{align*}
where the first estimate is obtained 
 from the fact that $\int G^2(m_s'(\cdot),j,t)w(t)dt=O(b_n/(nb_n^3))$. 
This  leads to the estimate
\begin{align}
{\rm Var} (nb_n^{9/2}D_1)=O\Big(\frac{1}{nb_n}\Big).\label{eq39}
\end{align}
For the term  $D_3$ in the decomposition \eqref{decom2} it follows that
\begin{align*}
\E(D_3^2)&=4\sum_{1\leq i\leq n}\Big(\int \tilde G(m_1'(\cdot),i,t)\tilde G(m_2'(\cdot),i,t)w(t)dt\Big)^2 \notag\\
&=\frac{4(\int vK'_d(v)dv)^4}{n^4b^{12}}\sum_i\Big(\int (((m_1')^{-1})')^2(t)(((m_2')^{-1})'(t) (K^\circ)'\Big(\frac{i/n-(m'_1)^{-1}(t)}{b_n}\Big)\notag\\& (K^\circ)'\Big(\frac{i/n-(m'_2)^{-1}(t)}{b_n}\Big)w(t)dt\Big)^2\sigma^2_1(i/n)\sigma^2_2(i/n)=O((n^3b_n^{11})^{-1})
\end{align*}
Hence, 
\begin{align}
nb_n^{9/2}D_3=O_p\Big(\Big(\frac{1}{nb_n^2}\Big)^{1/2}\Big).\label{eq41}
\end{align}
Finally we investigate the term $D_2$ using  a central limit theorem for quadratic forms [see  \cite{de1987central}].  For this purpose define  the terms 
(note that $(K^\circ)'(\cdot)$ is symmetric  and has  bounded support)
\begin{align*}
V_{s,n}&=\sum_{1\leq i\neq j\leq n}\Big( (K^\circ)'\Big(\frac{i/n-(m_s')^{-1}(t)}{b_n}\Big)(K^\circ)'\Big(\frac{j/n-(m_s')^{-1}(t)}{b_n}\Big)\sigma_s(
\frac{i}{n})\sigma_s(\frac{j}{n})(((m_s')^{-1})'(t))^4w(t)dt\Big)^2\\
&=n^2\int_0^1\int_0^1\Big(\int_{\mathbb R} (K^\circ)'\Big(\frac{u-(m_s')^{-1}(t)}{b_n}\Big)(K^\circ)'\Big(\frac{v-(m_s')^{-1}(t)}{b_n}\Big)\sigma_s(
u)\sigma_s(v)(((m_s')^{-1})'(t))^4w(t)dt\Big)^2
\\&\times dudv(1+o(1))\\
&=n^2b_n^2\int_0^1\int_0^1\Big(\int_{\mathbb R} (K^\circ)'(y)(K^\circ)'(\frac{v-u}{b_n}+y)\sigma_s^2(u)w(m_s'(u))
(m_s''(u))^{-3}dy\Big)^2dudv(1+o(1))
\\&=n^2b_n^3\int ((K^\circ)'*(K^\circ)'(z))^2dz
\int(\sigma_s^2(u)w(m_s'(u))(m''_s(u))^{-3} )^2  du(1+o(1))~, 
\end{align*}
 then $\lim_{n\to \infty}  V_{s,n} / (n^2 b_n^3) $ exists $(s=1,2$) and  
\begin{align}
\label{eq43}
\lim_{n\rightarrow \infty}\frac{2(\int vK'_d(v)dv)^4n^2b_n^9}{(nb_n^3)^4}(V_{1,n}+V_{2,n})
= {V_T},
\end{align}
where the asymptotic variance $V_T$ is defined in Theorem \ref{thm32}.
Now  similar arguments as used in the proof of Lemma 4 in 
\cite{zhou2010nonparametric} show  that 
\begin{align}\label{eq42}
nb^{9/2}_nD_2\Rightarrow {N(0,V_{T}),} 
\end{align}
Combining this statement with \eqref{equiv}, \eqref{decom2}, \eqref{eq37}, \eqref{eq39}, and \eqref{eq41} 
finally gives
\begin{align}\label{new.35}
nb_n^{9/2}\int U_n^2(t)w(t)dt- B_n(0) 
\Rightarrow N(0, V_T).
\end{align}

\paragraph{Asymptotic properties of \eqref{Un(t)1}:}
Note  that 
\begin{align*}
\int U_n(t)((m_1^{-1})'-(m_2^{-1})')w(t)dt=\int (U_{n,1}(t)-U_{n,2}(t))((m_1^{-1})'-(m_2^{-1})')w(t)dt,
\end{align*}
where 
\begin{align*}
\int (U_{n,s}(t)((m_1^{-1})'-(m_2^{-1})')w(t)dt=\sum_{j=1}^nV_{j,s}\int G(m_s'(\cdot),j,t)(\rho_ng(t)+o(\rho_n))w(t)dt\notag\\
=O_p\Big(\Big(\frac{nb_n^2\rho_n^2}{n^2b^6_n}\Big)^{1/2}\Big)=O_p\Big(\frac{\rho_n}{(nb_n^4)^{1/2}}\Big).
\end{align*}
Observing that 
$$
\int G(m_s'(\cdot),j,t)\rho_ng(t)w(t)dt=O(\rho_nb_n/(nb_n^3)),
$$
the bandwidth conditions and the definition of $\rho_n$  give for $s=1,2$,
\begin{align}\label{new.38}
nb_n^{9/2}\int (U_{n,s}(t)((m_1^{-1})'-(m_2^{-1})')w(t)dt=O_p(b_n^{1/4}).
\end{align}

\paragraph{Asymptotic properties of \eqref{Un(t)2}:}
Note that  it follows  for the term \eqref{Un(t)2}
\begin{align}
\Big|\int U_{n,s}(t)R^\dag_n(t)w(t)dt\Big|\leq  \sup_t|R^\dag_n(t)|\int\sup_t\Big|\sum_{j=1}^nV_{j,s}G(m_s'(\cdot),j,t)\Big|w(t)dt\notag.
\end{align}
Observing  that $\displaystyle\sum_jG^2(m_s'(\cdot),j,t)=O(nb_n/(nb_n^3)^2)$
 we have
\begin{align*}
\sup_t\Big|\sum_{j=1}^nV_{j,s}G(m_s'(\cdot),j,t)\Big|=O_p\Big(\frac{\log^{1/2} n}{n^{1/2}b_n^{5/2}}\Big),
\end{align*}
and   the conditions on the  bandwidths and \eqref{new.28} yield 
\begin{align}\label{new.40}
nb_n^{9/2}&\Big|\int (U_{n,s}(t)(R^\dag_n(t))w(t)dt\Big|\notag
\\&=O_p\Big(\frac{\log^{1/2} n}{n^{1/2}b_n^{5/2}}\Big(\frac{\pi_n'}{h_d}+\frac{\pi_n^3}{h_d^2}+h_d+\frac{1}{Nh_d}\Big)nb_n^{9/2}\Big)=o_p(1).
\end{align}

\noindent
 The proof of  assertion \eqref{s1} is now completed using the decomposition \eqref{2019-32} and  the results  \eqref{a1},  \eqref{a2},  \eqref{a3}, 
 \eqref{new.35},
 \eqref{new.38} and \eqref{new.40}. 
 
 \subsubsection{Proof of \eqref{s2}}
 From the proof of \eqref{s1}  we have the  decomposition
 \begin{align*}
T_n-\tilde T_n&=\int (I_1(t)-I_2(t)+II(t))^2(\hat w(t)-w(t))dt
\\&=\int \left(U_n(t)+((m'_1)^{-1}(t))'-((m'_2)^{-1}(t))'+R_n^\dag(t)\right)^2(\hat w(t)-w(t))dt, 
\end{align*}
where quantities $I_s$, $II$,
$U_n(t)$ and $R_n^\dag(t)$ are defined in \eqref{is}, \eqref{ii}, and \eqref{new.620}. By the proof of \eqref{s1}, it then suffices to show that
\begin{align*}
    nb_n^{9/2}\int (U_n(t))^2(\hat w(t)-w(t))dt=o_p(1).
\end{align*}
Using the same arguments  as given in the proof of \eqref{s1}, this assertion follows from 
 \begin{align*}
    nb_n^{9/2}\int (\tilde U_n(t))^2(\hat w(t)-w(t))dt=o_p(1).
\end{align*}
where $\tilde U_n(t)$ is defined in \eqref{tildeUn(t)}.
Recalling the definition of  $a,b$  in  \eqref{aandb} 
it then follows (using similar arguments as given for the derivation of  \eqref{Jul1-11}) that
\begin{align*}\label{SuptildeU_n(t)}
    \sup_{t\in [a,b]}|\tilde U_n(t)|=O_p\left(\frac{\log n}{\sqrt{nb_{n}}b^2_n}\right).
\end{align*}
Furthermore, together with part (iii) of Proposition \ref{PropM} it follows that 
\begin{align}
    \int (\tilde U_n(t))^2(\hat w(t)-w(t))dt\leq \sup_{t\in [a,b]}|\tilde U_n(t)|)^2\int |\hat w(t)-w(t)|dt=O_p\left(\frac{\bar \omega_n\log^2 n}{nb_n^5}\right),
\end{align}
where $\bar \omega_n$   is defined in \eqref{omega}.
 Thus by our choices of bandwidth $nb_n^{9/2}\frac{\bar \omega_n\log^2 n}{nb_n^5}=o(1)$, from which result (ii) follows.
 
 \medskip
Finally, the assertion of the Theorem  \ref{thm32} follows from \eqref{s1}  and \eqref{s2}.\hfill $\Box$

\bigskip
\medskip

\noindent 	
{\bf Acknowledgments.}
This work has been supported in part by the Collaborative Research Center ``Statistical modelling of nonlinear
dynamic processes'' (SFB 823, Teilprojekt  A1, C1) of the German Research Foundation (DFG).  The authors are grateful  to Martina Stein, who typed parts of this paper with considerable technical expertise. 
\begin{small}
	\setlength{\bibsep}{2pt}

\end{small}
\end{document}